\documentclass[11pt]{amsart}

\usepackage{amsmath,amsthm,verbatim,amssymb,amsfonts,amscd, graphicx}
\usepackage{graphics}
\usepackage{soul}
\usepackage{tikz-cd}
\usepackage{pinlabel}
\usepackage{xcolor}
\usepackage{stackrel}
\usepackage{hyperref}%
\usepackage[margin=1.5in]{geometry}
\newtheorem{theorem}{Theorem}[section]
\newtheorem{corollary}[theorem]{Corollary}
\newtheorem{lemma}[theorem]{Lemma}
\newtheorem{proposition}[theorem]{Proposition}

\newtheorem*{NOLP}{Theorem \ref{thm:NOLP}}
\newtheorem*{waldhausen}{Theorem \ref{thm:waldhausen}}
 
\newtheorem{question}[theorem]{Question} 
 
\theoremstyle{definition}
\newtheorem{definition}[theorem]{Definition}
\newtheorem{remark}[theorem]{Remark}
\newtheorem{notation}[theorem]{Notation}
\newtheorem*{proposition*}{Proposition}

\theoremstyle{definition}
\newtheorem{algorithm}[theorem]{Algorithm}
\newtheorem{example}[theorem]{Example}
\newtheorem{exercise}[theorem]{Exercise}

\newcommand{\leftrarrows}{\mathrel{\raise.75ex\hbox{\oalign{%
  $\scriptstyle\leftarrow$\cr
  \vrule width0pt height.5ex$\hfil\scriptstyle\relbar$\cr}}}}
\newcommand{\lrightarrows}{\mathrel{\raise.75ex\hbox{\oalign{%
  $\scriptstyle\relbar$\hfil\cr
  $\scriptstyle\vrule width0pt height.5ex\smash\rightarrow$\cr}}}}
\newcommand{\Rrelbar}{\mathrel{\raise.75ex\hbox{\oalign{%
  $\scriptstyle\relbar$\cr
  \vrule width0pt height.5ex$\scriptstyle\relbar$}}}}

\makeatletter
\def\leftrightarrowsfill@{\arrowfill@\leftrarrows\Rrelbar\lrightarrows}
\newcommand{\xleftrightarrows}[2][]{\ext@arrow 3399\leftrightarrowsfill@{#1}{#2}}
\makeatother

\newcommand{\maggie}[1]{\marginpar{\textcolor{blue}{\tiny{#1}}}}

\newcommand{\red}[1]{\textcolor{red}{#1}}

\definecolor{violet}{rgb}{.6,.6,0}
\definecolor{green}{rgb}{.0,.8,0}

\newcommand{\T}{\mathcal{T}}
\newcommand{\D}{\mathcal{D}}

\newcommand{\K}{\mathcal{K}}
\let\int\relax
\newcommand{\int}{\mathring}

\newcommand{\boundary}{\partial}
\newcommand{\boundaryin}{\partial_{\kern0.05em \text{in}}}
\newcommand{\boundaryout}{\partial_{\text{out}}}

\DeclareMathOperator{\id}{{id}}

\DeclareMathOperator{\pt}{{pt}}

\newcommand\ttimes{\mathbin{%
    \stackrel{\sim}{\smash{\times}\rule{0pt}{0.7ex}}%
    }}

\title[Trisections of non-orientable 4-manifolds]{Trisections of non-orientable 4-manifolds}

\author{Maggie Miller}
\address{Department of Mathematics\\Massachusetts Institute of Technology\\  Cambridge, MA 02142, USA}
\urladdr{https://math.mit.edu/~maggiehm/}
\email{maggiehm@mit.edu}

\author{Patrick Naylor}
\address{Department of Pure Mathematics\\University of Waterloo\\Waterloo, ON N2L 3G1, Canada}
\urladdr{https://patricknaylor.org}
\email{patrick.naylor@uwaterloo.ca}

\thanks{MM is supported by NSF Grant DMS-2001675. Earlier in this project, she was supported by NSF Grant DGE-1656466 at Princeton University. PN is supported by an NSERC CGS-D scholarship.}

\begin{document}
 
\maketitle

\begin{abstract}
We study trisections of smooth, compact non-orientable 4-manifolds, and introduce trisections of non-orientable 4-manifolds with boundary. In particular, we prove a non-orientable analogue of a classical theorem of Laudenbach-Po\'enaru. As a consequence, trisection diagrams and Kirby diagrams of closed non-orientable 4-manifolds exist. We discuss how the theory of trisections may be adapted to the setting of non-orientable 4-manifolds with many examples. 
\end{abstract}

\section{Introduction}

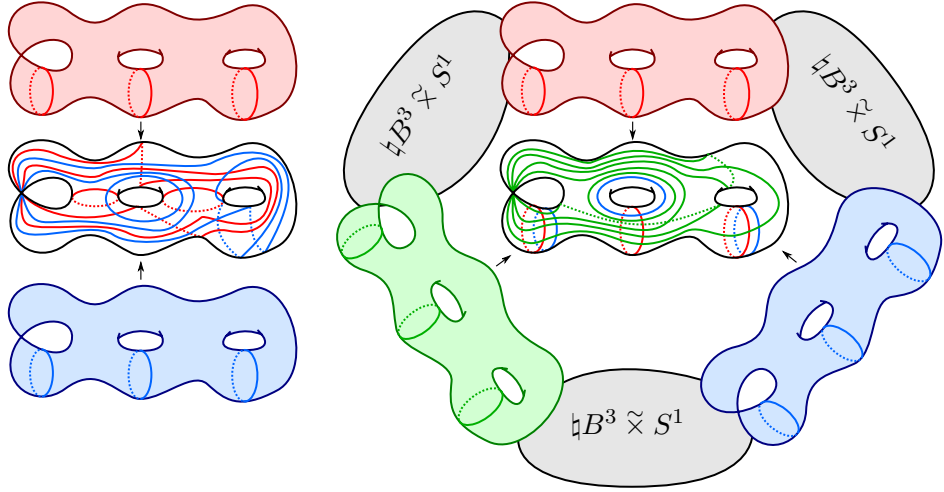
\begin{figure}
    \centering
\begingroup%
  \makeatletter%
  \providecommand\color[2][]{%
    \errmessage{(Inkscape) Color is used for the text in Inkscape, but the package 'color.sty' is not loaded}%
    \renewcommand\color[2][]{}%
  }%
  \providecommand\transparent[1]{%
    \errmessage{(Inkscape) Transparency is used (non-zero) for the text in Inkscape, but the package 'transparent.sty' is not loaded}%
    \renewcommand\transparent[1]{}%
  }%
  \providecommand\rotatebox[2]{#2}%
  \newcommand*\fsize{\dimexpr\f@size pt\relax}%
  \newcommand*\lineheight[1]{\fontsize{\fsize}{#1\fsize}\selectfont}%
  \ifx\svgwidth\undefined%
    \setlength{\unitlength}{352.12639275bp}%
    \ifx\svgscale\undefined%
      \relax%
    \else%
      \setlength{\unitlength}{\unitlength * \real{\svgscale}}%
    \fi%
  \else%
    \setlength{\unitlength}{\svgwidth}%
  \fi%
  \global\let\svgwidth\undefined%
  \global\let\svgscale\undefined%
  \makeatother%
  \begin{picture}(1,0.52094518)%
    \lineheight{1}%
    \setlength\tabcolsep{0pt}%
    \put(0,0){\includegraphics[width=\unitlength,page=1]{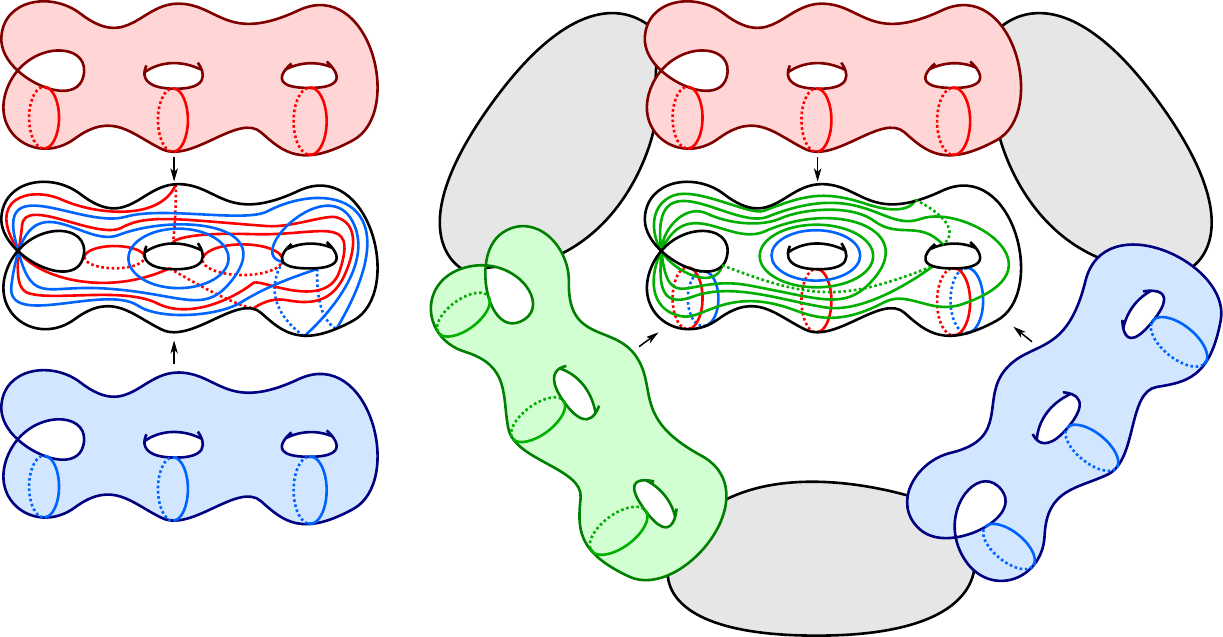}}%
    \put(0.86088376,0.45341934){\color[rgb]{0,0,0}\rotatebox{-54.03255599}{\makebox(0,0)[lt]{\lineheight{1.25}\smash{\begin{tabular}[t]{l}$\natural B^3\ttimes S^1$\end{tabular}}}}}%
    \put(0.41478606,0.35473719){\color[rgb]{0,0,0}\rotatebox{54.4237865}{\makebox(0,0)[lt]{\lineheight{1.25}\smash{\begin{tabular}[t]{l}$\natural B^3\ttimes S^1$\end{tabular}}}}}%
    \put(0.60209432,0.05763245){\color[rgb]{0,0,0}\makebox(0,0)[lt]{\lineheight{1.25}\smash{\begin{tabular}[t]{l}$\natural B^3\ttimes S^1$\end{tabular}}}}%
  \end{picture}%
\endgroup%

    \caption{{\bf{Left}}: a schematic of a Heegaard splitting of a (non-orientable) 3-manifold. Two (non-orientable) handlebodies are glued along their boundaries. The gluing map is specified by the image of curves in each handlebody that bound disks. Therefore, the 3-manifold is determined by a surface containing two sets of curves. {\bf{Right}}: a schematic of a trisection of a (non-orientable) 4-manifold. Three (non-orientable) 3-dimensional handlebodies are glued along their boundaries. The gluing maps are determined by the image of curves in each handlebody that bound disks. By \cite{laupoe72} or Theorem \ref{thm:NOLP}, the rest of the 4-manifold can be filled in uniquely up to diffeomorphism. Therefore, the 4-manifold is determined by a surface containing three sets of curves.}
    \label{fig:introduction}
\end{figure}

In this paper, we study trisections of smooth, compact, non-orientable 4-manifolds (possibly with boundary). A trisection is a decomposition of a 4-manifold into three 4-dimensional handlebodies, whose triple intersection is a properly embedded surface. This is analogous to a Heegaard splitting of a 3-manifold, which is a decomposition into two 3-dimensional handlebodies whose intersection is a properly embedded surface. A schematic is given in Figure \ref{fig:introduction}.

Trisections of orientable 4-manifolds were introduced by Gay and Kirby \cite{gaykirby}. Trisections of closed, non-orientable 4-manifolds were later studied by Rubinstein and Tillmann \cite{multisections} and by Spreer and Tillmann \cite{spretil20}. We will show that many theorems about trisections of orientable 4-manifolds also hold for non-orientable 4-manifolds. In particular, we prove the following two theorems, which may be of independent interest.

\begin{NOLP}
Let $h:\#_p S^2\ttimes S^1\to \#_p S^2\ttimes S^1$ be a diffeomorphism. Then there is a diffeomorphism $H:\natural_pB^3\ttimes S^1\to \natural_p B^3\ttimes S^1$ such that $H\vert_{\partial Y_p}=h$.
\end{NOLP}

\begin{waldhausen}
Fix $g\ge k\ge 0$. Any genus $g$ Heegaard surface of $\#_k S^2\ttimes S^1$ is equivalent to the result of stabilizing the standard Heegaard surface $g-k$ times.
\end{waldhausen}

See Section \ref{sec:background} for detailed exposition on trisections of closed 4-manifolds and Section \ref{sec:relative} for 4-manifolds with boundary. As an introduction, we will state some notable properties of trisections:
\begin{enumerate}
    \item Any smooth, compact 4-manifold $X^4$ admits a trisection.
    \item Any two trisections of a closed 4-manifold become isotopic after connect-summing the trisections with some number of genus one trisections of $S^4$. (This is called {\emph{stabilization}} or {\emph{interior stabilization}}.)
    \item  Any two relative trisections of a compact, orientable 4-manifold $X^4$ with $\boundary X^4\neq\emptyset$ are related by a finite sequence of interior stabilizations and two additional moves. It is unclear how to modify this statement to hold for non-orientable 4-manifolds.
    \item\label{diagramitem} A trisection can be described by a {\emph{trisection diagram}} consisting of a surface $\Sigma$ and three sets of simple closed curves on $\Sigma$, usually denoted $\alpha$, $\beta$, and $\gamma$. This diagram is determined by the trisection up to automorphism of $\Sigma$ and slides of $\alpha$, $\beta$, and $\gamma$.
    \item\label{surfaceitem} Any surface smoothly embedded in a trisected 4-manifold can be put into {\emph{bridge position}} with respect to the trisection (similarly to how links in a Heegaard split 3-manifold can be put into bridge position). This allows one to encode a surface $S$ in a 4-manifold via a trisection diagram with additional dots and arcs describing the intersection of $S$ with various pieces of the trisection. Moreover, there is a perturbation operation on a surface in bridge position so that if two isotopic surfaces $S_1$, $S_2$, are in bridge position, then after a sequence of perturbations, $S_1$ and $S_2$ are isotopic by an isotopy that fixes the pieces of the trisection setwise. (That is, bridge position is unique up to perturbation and deperturbation, just as in a 3-manifold.)
    \item\label{coveritem} If $\D$ is a trisection diagram of a 4-manifold $X^4$, then regular covers of $\D$ are trisection diagrams of covers of $X^4$. Similarly, trisection diagrams indicating a surface $S$ in bridge position may be modified to produce diagrams of covers branched over $S$.
\end{enumerate}

Items \ref{diagramitem}, \ref{surfaceitem}, and \ref{coveritem} give compelling reasons to study trisections of 4-manifolds. A trisection diagram is similar to a Kirby diagram of a 4-manifold, but contains an additional $\mathbb{Z}/3$ symmetry not found in Kirby diagrams: one can turn a Kirby diagram upside-down, but one can ``rotate a trisection through $2\pi/3$'' by exchanging the names of the three pieces of the trisection. Diagrams of surfaces in 4-manifolds via trisections are extremely combinatorial and perhaps more acccessible for computation than standard banded unlinks in Kirby diagrams \cite{hugkimmil20}, and taking covers or branched covers from a trisection is significantly easier than performing the same operation on a Kirby diagram.

With at least some understanding of a trisection diagram, we can explain the importance of Theorems \ref{thm:NOLP} and \ref{thm:waldhausen}.
\begin{itemize}
    \item Theorem \ref{thm:NOLP} allows one to construct a 4-manifold from a non-orientable trisection diagram. Additionally, Theorem \ref{thm:NOLP} ensures that Kirby diagrams exist for non-orientable 4-manifolds (without having to indicate the attaching spheres of 3-handles).
    \item Theorem \ref{thm:waldhausen} ensures that a closed, trisected 4-manifold determines a trisection diagram up to surface automorphism and slides of curves.
\end{itemize}

\subsection*{Organization}
The sections of this paper are as follows.
\begin{enumerate}
    \item[\S\ref{sec:background}:] We review some background on trisections of closed 4-manifolds.
    \item[\S\ref{sec:LP}:] We prove Theorem \ref{thm:NOLP}. We use this to define trisection diagrams for non-orientable 4-manifolds, and give some applications to diffeomorphisms of $\#_p S^2\ttimes S^1$. As a consequence, we conclude that Kirby diagrams for closed non-orientable 4-manifolds exist.
    \item[\S\ref{sec:waldhausen}:] We prove Theorem \ref{thm:waldhausen} and conclude that trisections determine trisection diagrams.
    \item[\S\ref{sec:nonorientabletrisections}:] We produce some non-orientable trisections and show how one may take orientation double covers of trisected 4-manifolds.
    \item[\S\ref{sec:relative}:] We define relative trisections of 4-manifolds with boundary, which is a considerably more involved definition than in the closed case. We discuss some commonly performed operations on relative trisections: gluing two relative trisections, converting a relative trisection diagram to a Kirby diagram, and the reverse.
    \item[\S\ref{sec:bridge}:] We explain why surfaces smoothly embedded in a trisected non-orientable manifold can be put into bridge position, and why any two isotopic surfaces in bridge position with respect to a non-orientable trisection are related by a sequence of perturbations and deperturbations. We give a few examples of surfaces in bridge position.
\end{enumerate}

\subsection*{Acknowledgements}
We thank Fran\c{c}ois Laudenbach for extremely helpful correspondence, as mentioned briefly in Section \ref{sec:LP}. The second author would like to thank his graduate advisor, Doug Park. This project began when the first author visited the University of Waterloo in December, 2019.

\section{Background: trisections of 4-manifolds}\label{sec:background}

In this section, we briefly recall the definition of a trisection of a 4-manifold, first introduced by Gay and Kirby in \cite{gaykirby}.

\begin{definition}
An \emph{handlebody} of genus $g$ is a compact manifold admitting a handle decomposition with a single 0-handle and $g$ 1-handles. Handlebodies will generally \emph{not} be assumed to be orientable. 
\end{definition}

\begin{definition}
Suppose that $X^4$ is a smooth, closed, and connected 4-manifold. A $(g;k)$-\emph{trisection} of $X$ is a decomposition $X=X_1\cup X_2\cup X_3$ such that:
\begin{itemize}
    \item $X_i$ is diffeomorphic to a 4-dimensional handlebody of genus $k$;
    \item $X_i\cap X_j$ is diffeomorphic to a 3-dimensional handlebody of genus $g$;
    \item $X_1\cap X_2\cap X_3\cong \Sigma_g$, a closed surface of genus $g$.
\end{itemize}
\end{definition}

This definition does not require $X$ to be orientable. We often refer to the triple intersection $\Sigma$ as the \emph{central surface} and the 4-dimensional handlebodies as \emph{sectors} of the trisection $(X_1,X_2,X_3)$. Note that the central surface induces a genus $g$ Heegaard splitting of the boundary of each sector. By a Heegaard splitting of a non-orientable 3-manifold, we simply mean a decomposition into two (necessarily non-orientable) handlebodies which intersect in a common surface. 
By \cite[Proposition 5]{spretil20}, all handlebodies (and central surface) of a trisection of $X$ are orientable if and only if $X$ is orientable. In particular, if $X$ is non-orientable, then none of the sectors can be 4-balls. Our notation will differ slightly from \cite{spretil20}: here, the \emph{genus} of a trisection is the the genus of the 3-dimensional handlebodies bounded by the central surface. When $X$ is orientable, we are assured of existence and uniqueness of trisections of $X$ up to a natural stabilization operation by \cite{gaykirby}. For any closed 4-manifold, Rubinstein and Tillmann \cite{multisections} give a proof of existence and stable equivalence that does not require orientability. 

There is also a notion of a {\emph{relative trisection}} of a compact 4-manifold with boundary. A relative trisection is similarly unique up to interior and relative stabilization operations \cite{gaykirby,castro,CIMT}. We discuss this further in Section \ref{sec:relative}.

One of the main appeals of trisections is that they provide a diagrammatic calculus with which one can study 4-manifolds.

\begin{definition}
A $(g;k)$- \emph{trisection diagram} is a tuple $(\Sigma;\alpha,\beta,\gamma)$, where $\Sigma$ is a closed orientable surface of genus $g$, and $\alpha$, $\beta$, and $\gamma$ are collections of $g$ embedded curves on $\Sigma$ such that:
\begin{itemize}
    \item Each of $\alpha$, $\beta$, and $\gamma$ is a \emph{cut system of curves} for $\Sigma$;
    \item Each pair of curves is \emph{standard}, i.e.\  each of $(\Sigma;\alpha,\beta)$, $(\Sigma;\beta,\gamma)$, and $(\Sigma;\gamma,\alpha)$ is a Heegaard diagram for $\#_kS^2\times S^1$. 
\end{itemize}
\end{definition}

\noindent A trisection diagram determines a trisection of a closed 4-manifold via a theorem of Laudenbach-Po\'enaru \cite{laupoe72} in the usual way; see \cite{gaykirby}.

\section{A non-orientable version of Laudenbach—Po\'enaru}\label{sec:LP}

In this section, we give an analogue of a theorem of Laudenbach and Po\'enaru \cite{laupoe72} for non-orientable 4-dimensional handlebodies, which may be of independent interest. While we suspect that this theorem might be known to some experts, it does not appear in the literature and is critical for diagrammatic descriptions of closed, non-orientable 4-manifolds.

\begin{theorem}\label{thm:NOLP}
Let $h:\#_p S^2\ttimes S^1\to \#_p S^2\ttimes S^1$ be a diffeomorphism. Then there is a diffeomorphism $H:\natural_pB^3\ttimes S^1\to \natural_p B^3\ttimes S^1$ such that $H\vert_{\partial Y_p}=h$.
\end{theorem}

The proof will use the same strategy as \cite{laupoe72}. The first step is to reduce to the case that $h$ acts homotopically trivially. After this, one uses Laudenbach's theorem on homotopy and isotopy for 2-spheres in 3-manifolds \cite{lau73} together with Cerf's theorem \cite{cer68} to isotope $h$ to a diffeomorphism which clearly extends to $\natural^p B^3\ttimes S^1$. In particular, we make use of the following theorem of Laudenbach \cite{lau73}\footnote{We are grateful to Fran\c{c}ois Laudenbach for corresponding with us about this remarkable theorem.}

\begin{theorem}[\cite{lau73}]\label{laudenbachspheretheorem}
Let $S$ and $T$ be embedded 2-spheres in a 3-manifold $V$, which may be non-orientable and have nonempty boundary. Then if $S$ and $T$ are homotopic, they are isotopic.
\end{theorem}

We will denote $V_p=\#_pS^2\ttimes S^1$, and $Y_p=\natural_p B^3\ttimes S^1$. Observe that we have the following commutative triangle of natural homomorphisms:

\vspace{2mm}
\begin{center}
\begin{tikzcd}
 & \pi_0(\textnormal{Diff}(Y_p)) \arrow[dl,"A"'] \arrow[dr,"B"] & \\
 \textnormal{Aut}(\pi_1(Y_p)) \arrow[rr,"\textnormal{id}"'] &  &  \textnormal{Aut}(\pi_1(V_p))\\
\end{tikzcd}
\end{center}
\vspace{-5mm}

\begin{lemma}
The maps $A$ and $B$ are surjective. 
\end{lemma}

\begin{proof}
It is sufficient to prove that $A$ is surjective. The surjectivity of $A$ follows from the fact that $Y_p$ admits a handle decomposition with one $0$-handle and $p$ $1$-handles. That is, $\pi_1(Y_p))$ is free with generators $a_1,\ldots, a_p$ corresponding to cores of these $1$-handles, in order. The group $\textnormal{Aut}(\pi_1(Y_p))$ is generated by automorphisms of the form:
\begin{itemize}
    \item[(a)] $a_i\mapsto a_i^{-1}$ for some $i$, i.e.\  inverting a generator;
    \item[(b)] $a_i\mapsto a_j$, $a_j\mapsto a_i$ for some $i\neq j$, i.e.\  interchanging two generators;
    \item[(c)] $a_i\mapsto a_ia_j$ for some $i\neq j$, i.e.\  multiplying a generator on the right by another generator. 
\end{itemize}

The first two automorphisms are realizable by diffeomorphisms of $Y^p$ in which the feet of $1$-handles exchange positions on the boundary of the $0$-handle. The third automorphism is induced by the diffeomorphism sliding the $i$-th 1-handle over the $k$-th $1$-handle.
\end{proof}

By composing with a diffeomorphism of $Y_p$, we may now assume that $h^\#_1:\pi_1(V_p)\to \pi_1(V_p)$ is the identity. 

\begin{lemma}\label{idonpi2}
Suppose that a diffeomorphism $h:V_p\to V_p$ is such that $h^\#_1:\pi_1(V_p)\to \pi_1(V_p)$ is the identity. Then $h^\#_2:\pi_2(V_p)\to \pi_2(V_p)$ is also the identity. 
\end{lemma}

\begin{proof}
Since $h_1^\#$ is the identity map, $h$ lifts to a diffeomorphism $\tilde{h}:\#_{2p-1}S^2\times S^1\to \#_{2p-1}S^2\times S^1$ of the orientation double cover of $Y_p$. Then $\tilde{h}_1^\#$ is either the identity or an involution. By composing with the deck transformation if necessary, we may assume that $\tilde{h}$ is orientation preserving, so $\tilde{h}_1^\#$ is the identity. By Lemma 3 of \cite{laupoe72}, $\tilde{h}^\#_2$ is the identity. Since the covering map induces the identity map on $\pi_2$, we conclude that $h^\#_2$ is the identity map.
\end{proof}

Now, let $S_1,\dots,S_p \subset V_p$ denote $p$ disjoint, homologically independent 2-spheres of the form $S^2\times \{\text{pt}\}$. By Lemma \ref{idonpi2}, we conclude that $h(S_i)$ and $S_i$ are homotopic. We will show that the systems of spheres $S_1,\ldots, S_p$ and $h(S_1),\ldots, h(S_p)$ are isotopic. The following proposition will be useful.

\begin{proposition}\label{injectiveprop}
Let $M$ be a (possible non-orientable) 3-manifold whose boundary is a disjoint union of 2-spheres, and let $N$ be a 3-manifold obtained by gluing $S^2\times I$ to $M$ along two boundary components $B_1,B_2$ of $M$. 
Suppose $S_1$ and $S_2$ are 2-spheres in $M$ that are not homotopic to each other or to $B_1$ or $B_2$. Then $S_1$ and $S_2$ are also not homotopic in $N$.
\end{proposition}
\begin{proof}





Let $\widetilde{M}$ denote the universal cover of $M$, and $\widetilde{B}_i$ the set of boundary components of $\widetilde{M}$ that are lifts of $B_i$. Then the universal cover $\widetilde{N}$ of $N$ is obtained from $\widetilde{M}$ by gluing a different copy of $S^2\times I$ along one boundary to each element of $\widetilde{B}_1,\widetilde{B}_2$, then gluing copies of $\widetilde{M}$ to each remaining boundary of an $S^2\times I$ (along a sphere in $\widetilde{B}_1$ or $\widetilde{B}_2$, appropriately), and iterating infinitely. This construction yields an inclusion of $\widetilde{M}$ in $\widetilde{N}$.

Suppose $S_1$ and $S_2$ are homotopic in $N$, so that they have lifts $L_1$ and $L_2$ in $\widetilde{N}$ that cobound an immersed copy $W$ of $S^2\times I$ in $\widetilde{N}$. 


Let $W'$ be comprised of the one or two (closures of) components of $W\setminus(\widetilde{B}_1\cup\widetilde{B}_2)$ with $L_1$ or $L_2$ as a boundary component. Since $S_1$ and $S_2$ are not homotopic in $M$, $W'$ is not all of $W$.

Suppose $W\setminus W'$ meets some copy of $\widetilde{M}$ in $\widetilde{N}$. Then based spheres representing $B_1,B_2$ in $M$ are not free generators of $\pi_2(M)$, implying $M\cong S^2\times I$, contradicting the existence of non-homotopic spheres $S_1,S_2$. Therefore, $W\setminus W'\cong S^2\times I$ and $W'$ has two components: an immersed copy of $S^2\times I$ from $L_1$ to a sphere in $\widetilde{B}_1\cup\widetilde{B}_2$ and a similar copy of $S^2\times I$ from a sphere in $\widetilde{B}_1\cup\widetilde{B}_2$ to $L_2$. This implies that $S_1$ and $S_2$ can each be homotoped in $M$ to one of $B_1$, $B_2$, which is a contradiction.


%

\end{proof}

\begin{lemma}\label{fixesSi}
Let $h:V_p\to V_p$ be a diffeomorphism. If $h(S_i)$ is homotopic to $S_i$ for all $i$, then $h$ is isotopic to a diffeomorphism which is the identity on each $S_i$. 
\end{lemma}

\begin{proof} 
Inductively, assume that $h$ fixes $S_j$ pointwise for all $j<k$. We will show that we can isotope $h$ in the complement of $S_1,\ldots, S_{k-1}$ to also fix $S_k$ pointwise.

Note that $[S_k]=[h(S_k)]$ is linearly independent from $[S_1],\ldots,[S_{k-1}]$ in $H_2(V_p)$. Then by repeated application of Proposition \ref{injectiveprop}, since $S_k$ and $h(S_k)$ are homotopic in $V_p$, there are homotopic in $V_p\setminus(S_1\sqcup\cdots\sqcup S_{k-1})$. Then by Theorem \ref{laudenbachspheretheorem}, $h$ can be isotoped in the complement of $S_1,\ldots, S_{k-1}$ to fix $S_k$ setwise and with $h|_{S_k}$ orientation-preserving, and then further to fix $S_k$ pointwise.




\end{proof}


For $\alpha\in \pi_1(SO(3))$, let $H_\alpha$ denote the diffeomorphism of $S^2\times I$ defined by $H_\alpha(x,t)=(\alpha(t)x,t)$. Since $\pi_1(SO(3))=\mathbb{Z}/2\mathbb{Z}$, there are only two possibilities for such a map up to isotopy, depending on the homotopy class of $\alpha$. If $S\subset V_p$ is an embedded 2-sphere, we will use $H_\alpha(S)$ to denote $H_\alpha$ applied to $\nu(S)\cong S^2\times I$. 

\begin{lemma}\label{lem:reducetotwists}
Let $h:V_p\to V_p$ be a diffeomorphism. If $h$ is the identity on each $S_i$, then $h$ is isotopic to a diffeomorphism of the form $H_{\alpha_{p}}(S_p)\circ\cdots\circ H_{\alpha_{1}}(S_1)$.
\end{lemma}
\begin{proof}
First, note that cutting $V_p$ along $S_1,\dots,S_p$ gives a $2p$-punctured 3-sphere $S(2p)$. Moreover, since $h$ is the identity on each $S_i$, $h$ induces an orientation preserving diffeomorphism $f:S(2p)\to S(2p)$, which is the identity on all boundary components. By Cerf's theorem \cite{cer68}, $f$ is isotopic to the identity map except possibly in a collar neighborhood of the boundary components, where $f$ may differ from the identity by some twist $H_\alpha$. Re-gluing the boundary components of the 3-sphere to form $V_p$, we conclude that $h$ is isotopic to a composition of twists along $S_1,\dots,S_p$ as above.
\end{proof}

Now, any diffeomorphism of $V_p$ of the form $H_\alpha(S_i)$ extends to $Y_p$, since the rotation of $S_i$ extends to the 3-ball that it bounds in $Y_p$. Thus $h$ extends to $Y_p$, and this completes the proof of Theorem~\ref{thm:NOLP}. \qed





\subsection{Trisection diagrams}

Using Theorem \ref{thm:NOLP}, we can define trisection diagrams for non-orientable 4-manifolds.

\begin{definition}
A \emph{$(g;k)$-trisection diagram} is a tuple $(\Sigma_g;\alpha,\beta,\gamma)$, where $\Sigma_g$ is a non-orientable surface of genus $g$, and $\alpha,\beta$, and $\gamma$ are three collections of $g$ curves such that:
\begin{itemize}
    \item Each of $\alpha,\beta$, and $\gamma$ is a \emph{cut system of curves} for $\Sigma_g$, that is, their complement in $\Sigma_g$ is planar, and each curve has an annular neighbourhood;
    \item Each of $(\Sigma_g;\alpha,\beta),(\Sigma_g;\beta,\gamma),(\Sigma_g;\gamma,\alpha)$ is a Heegaard diagram for $\#_kS^2\ttimes S^1$.
\end{itemize}
\end{definition}

Such a trisection diagram describes a trisected 4-manifold in the usual way. Let $V_\alpha,V_\beta,V_\gamma$ denote the 3-dimensional handlebodies prescribed by the $\alpha,\beta,\gamma$ curves respectively. Beginning with $\Sigma_g\times D^2$, attach $V_\alpha \times I, V_\beta\times I$ and $V_\gamma\times I$ to $\partial (\Sigma_g\times D^2)=\Sigma_g\times S^1$ along $\Sigma_g\times[-\epsilon,\epsilon],\Sigma_g\times[2\pi/3-\epsilon,2\pi/3+\epsilon],\Sigma_g\times[4\pi/3-\epsilon,4\pi/3+\epsilon]$, respectively. The resulting manifold has three boundary components, each diffeomorphic to $\#_k S^2\ttimes S^1$. By Theorem \ref{thm:NOLP}, we may uniquely (up to diffeomorphism) fill each component with $\natural^k B^3\ttimes S^1$ to get a closed 4-manifold.

A diagrammatic characterization for Heegaard splittings of $\#_kS^2\ttimes S^1$ will follow from the corresponding analogue of Waldhausen's theorem in Section~\ref{thm:waldhausen}.

\subsection{Other applications}
In this subsection, we discuss some other applications of Theorem \ref{thm:NOLP}.

\begin{theorem}\label{homotopicimpliesisotopicS1xS2}
Let $h:\#_p S^2\ttimes S^1\to \#_p S^2\ttimes S^1$ be a diffeomorphism homotopic to the identity. Then $h$ is isotopic to the identity.
\end{theorem}

\begin{proof}
This follows immediately from Lemma \ref{lem:reducetotwists}, and the following proposition due to Laudenbach. For the convenience of the reader (especially one who does not speak French), we give a proof below. 
\end{proof}

\begin{proposition}[\cite{laudenbach74}, Appendix II]\label{prop:laudenbachsprop}
Let $M$ be a (possibly non-orientable) 3-manifold, and let $S_1,\dots,S_p\subset M$ be embedded 2-spheres. Suppose that $S_1$ does not separate $M- (S_2\cup\cdots\cup S_p)$ and that the map $H:=H_{\alpha_1}(S_1)\circ\cdots \circ H_{\alpha_p}(S_p)$ is homotopic to the identity. Then $\alpha_1=0$. 
\end{proposition}

\begin{proof}[Proof of Proposition \ref{prop:laudenbachsprop}]

Since $S_1$ is non-separating in $M- (S_2\cup\cdots\cup S_p)$, there exists a circle $\gamma$ in $M- (S_2\cup\cdots\cup S_p)$ intersecting $S_1$ transversely once.

\vspace{2mm}
\noindent {\bf{Case 1: the normal bundle of $\gamma$ is orientable.}}
\vspace{2mm}


Let $T_0$ be a trivialization of the normal bundle of $\gamma$, and let $T_1=T_0\circ dH$ be the corresponding trivialization induced by $H$. Since $H$ is homotopic to the identity, there exists a framed annulus $A$ in $M\times\mathbb{R}$ with boundary $((\gamma,T_0)\times 0)\sqcup(-(\gamma,T_1)\times 1)$. The framing of $A$ extends $T_0$ and $T_1$ to trivializations $\widetilde{T}_0$ and $\widetilde{T}_1$ of the normal bundle $\nu_{M\times\mathbb{R}}(\gamma)$ of $\gamma$ in $M\times \mathbb{R}$, with $\widetilde{T}_0$ and $\widetilde{T}_1$ homotopic.

Suppose $\alpha_1=1$. From the construction of $T_0$ and $T_1$, we see that when we view $\widetilde{T}_0$ and $\widetilde{T}_1$ as elements of $\pi_1(SO(3))$, they differ by a $2\pi$ rotation about some axis. This is exactly the generator of $\pi_1(SO(3))\cong\mathbb{Z}/2$, contradicting $\widetilde{T}_1\sim\widetilde{T}_0$. Therefore, $\alpha_1=0$.

\vspace{2mm}
\noindent{\bf{Case 2: the normal bundle of $\gamma$ is non-orientable.}}
\vspace{2mm}

Let $\widetilde{M}\to M$ denote the orientation double cover of $M$ and let $\tau:\widetilde{M}\to\widetilde{M}$ the deck transformation. For $i=1,\ldots, p$, let $\widetilde{S}_i$ be one lift of $S_i$ to $\widetilde{M}$. Then $H$ lifts to the map \[\widetilde{H}=H_{\alpha_1}(\widetilde{S}_1)\circ\cdots\circ H_{\alpha_p}(\widetilde{S}_p)\circ H_{\alpha_1}(\tau\widetilde{S}_1)\circ\cdots\circ H_{\alpha_p}(\tau\widetilde{S}_p).\]

Since $H$ is homotopic to the identity, $\widetilde{H}$ is homotopic to either the identity or $\tau$. But $\tau$ is orientation-reversing while $\widetilde{H}$ is orientation-preserving, so $\widetilde{H}$ must be homotopic to the identity.

We would like to apply the proposition to $\widetilde{H}$, since we have already shown the claim to be true when the ambient 3-manifold is orientable. Unfortunately, we only know that $(\widetilde{S}_1\cup\tau\widetilde{S}_1)$ is nonseparating in the complement of $\widetilde{S}_2\cup\cdots\cup\widetilde{S}_p\cup\tau\widetilde{S}_2\cup\cdots\cup\tau\widetilde{S}_p$, which does not match the hypothesis. We must perform some surgery on $\widetilde{M}$ to arrange for $\widetilde{S}_1$ to be nonseparating in the complement of $\widetilde{S}_2\cup\cdots\cup\widetilde{S}_p\cup\tau\widetilde{S}_1\cup\tau\widetilde{S}_2\cup\cdots\cup\tau\widetilde{S}_p$. To this end, isotope $H$ to fix a small closed ball $B$ about a point $x_0$ pointwise. 

\begin{lemma}\label{lem:laudenbachappendixcase2}
There exists a homotopy $h_t$ from $H$ to $\id$ so that $h_t^{-1}(B)=B$ and $h_t^{-1}(x_0)=x_0$ for all $t$. (That is, $h_t$ fixes $B$ setwise and $x_0$ pointwise.)
\end{lemma}

\begin{proof}[Proof of Lemma \ref{lem:laudenbachappendixcase2}]
Let $g_t$ be a homotopy from $H$ to the identity. Then as $t$ ranges from $0$ to $1$, $g_t(x_0)$ traces out a loop $\eta$. If $[\eta]=0\in\pi_1(M,x_0)$, then we may obtain a new homotopy $h_t$ as desired by contracting $\eta$ to a point.

\vspace{2mm}
\noindent{\bf{Case 2(a): $M\cong S^2\ttimes S^1$.}}
\vspace{2mm}

Since $H$ is homotopic to the identity, $\eta$ lifts to a loop in $\widetilde{M}$. Therefore, $[\eta]\in2\mathbb{Z}\subset\mathbb{Z}=\pi_1(M)$.

Consider the homotopy $f_t$ of $S^2\ttimes S^1$ from the identity to the identity that sends each $\pt\times S^2$ twice around the $S^1$ direction as $t$ goes from $0$ to $1$. Then $f_t(x_0)$ traces out a curve in the homotopy class $2\in\mathbb{Z}=\pi_1(M)$. By concatenating $g_t$ with some power of $f_t$ or its inverse, we may arrange that $[\eta]=0$.

\vspace{2mm}
\noindent{\bf{Case 2(b): $M\not\cong S^2\ttimes S^1$.}}
\vspace{2mm}

Since $H$ is isotopic to the identity, the map on $\pi_1(M,x_0)$ given by conjugation with $[\eta]$ is trivial. That is, $[\eta]$ is in the center of $\pi_1(M)$. But since $M$ is a 3-manifold with $\pi_1(M)\not\cong\mathbb{Z}$, $\pi_1(M)$ has trivial center. Then $[\eta]=0$. This completes the proof of Lemma \ref{lem:laudenbachappendixcase2}.
\end{proof}

Thus, there exists a homotopy $h_t$ from $H$ to the identity fixing $x_0$ pointwise and $B$ setwise. Let $\widetilde{B}$ be a lift of $B$ to $\widetilde{M}$. Fix a homeomorphism $\phi:S^2\to \partial{\widetilde{B}}$ and consider the manifold
\[W:=\displaystyle\frac{\left(\widetilde{M}\setminus[\text{int}(\widetilde{B})\cup\text{int}(\tau\widetilde{B})]\right)\cup (S^2\times[0,1])}{(x\times 0\sim\phi(x)\text{ and } x\times 1\sim\tau\phi(x)\text{ for all $x\in S^2$})}.\]

The point of this construction is that $W\cong \widetilde{M}\# S^2\times S^1$, and we have arranged the surgery so that $\tilde{H}$ naturally extends to $W$ via the diffeomorphism $G:W\to W$ given by \[G(x)=\begin{cases}\widetilde{H}(x)&x\not\in S^2\times[0,1],\\x&x\in S^2\times[0,1].\end{cases}\]

Observe that following are true:
\begin{enumerate}
\item $\widetilde{S}_1,\ldots, \widetilde{S}_p,\tau \widetilde{S}_1,\ldots, \tau \widetilde{S}_p$ are 2-spheres embedded in $W$.
\item $\widetilde{S}_1$ is nonseparating in $W\setminus(\widetilde{S}_2\cup\cdots\cup \widetilde{S}_p\cup\tau \widetilde{S}_1\cup\tau \widetilde{S}_2\cup\cdots\cup\tau \widetilde{S}_p)$.
\item $G=H_{\alpha_1}(\widetilde{S}_1)\circ\cdots\circ H_{\alpha_p}(\widetilde{S}_p)\circ H_{\alpha_1}(\tau\widetilde{S}_1)\circ\cdots\circ H_{\alpha_p}(\tau\widetilde{S}_p).$
\item $G$ is homotopic to the identity (via $h_t$ on $\widetilde{M}\cap W$; here we use Lemma \ref{lem:laudenbachappendixcase2}).
\item $W$ is orientable.
\end{enumerate}

Then by applying Case 1 to the map $G$, we conclude that $\alpha_1=0$, completing the proof of Proposition \ref{prop:laudenbachsprop}.

\end{proof}

\begin{remark}\label{rem:diffeotopygroupS1xS2}
When $p=1$, Theorem \ref{homotopicimpliesisotopicS1xS2} gives another proof that the diffeotopy group of $S^2\ttimes S^1$ is equal to $\mathbb{Z}_2\oplus\mathbb{Z}_2$, generated by a reflection, and a twist about a 2-sphere fiber, which we will simply denote $\tau$. This was first computed by Kim and Raymond in \cite{kimray90}, who also gave the corresponding analogue of Theorem \ref{homotopicimpliesisotopicS1xS2}. One way to check that $\tau$ is \emph{not} homotopic to the identity is to first observe that $S^2\ttimes S^1=\partial( D^2\ttimes \mathbb{RP}^2)$, the non-trivial disk bundle over $\mathbb{RP}^2$ obtained as the complement of an orientation reversing loop in $\mathbb{RP}^4$. Doubling this bundle via either the identity or $\tau$ yields either $S^2\ttimes \mathbb{RP}^2$ or $\mathbb{RP}^4\#_{S^1}\mathbb{RP}^4$, the \emph{circle sum} of $\mathbb{RP}^4$ with itself. These are known to be homotopy inequivalent 4-manifolds by the complete classification of non-orientable 4-manifolds with fundamental group $\mathbb{Z}_2$ by Hambleton-Kreck-Teichner \cite{hamkretei94}, or the homotopy invariant for such manifolds given by Kim-Kojima-Raymond in \cite{kimkojray92}.
\end{remark}

In addition to the above results on Dehn twists on 2-spheres in 3-manifolds, Theorem \ref{thm:NOLP} allows one to use Kirby diagrams to describe closed, non-orientable 4-manifolds.
 
\begin{corollary}
Let $X^4$ be a closed, non-orientable 4-manifold. Fix a handle decomposition of $X$, and let $X_{(n)}$ denote the union of the $0$-, $1$-,\dots, $n$-handles of this decomposition. Then $X$ is determined up to diffeomorphism by $X_{(2)}$.
\end{corollary}
\begin{proof}
This follows immediately from Theorem \ref{thm:NOLP}, since $X\setminus X_{(2)}$ is a non-orientable $1$-handlebody.
\end{proof}
Thus, a diagram of $X_{(2)}$ along with the hypothesis that $X$ is closed determines $X$ up to diffeomorphism, as is the case for orientable 4-manifolds. Kirby diagrams have previously been used to study non-orientable 4-manifolds with boundary (we refer the reader to \cite{akbulutbook}). Akbulut has also used handle diagrams of closed, non-orientable 4-manifolds that include a description of the attaching spheres of 3-handles in e.g.\ \cite{akbulutfake}.

As is the case for orientable 4-manifolds \cite{kirbycalc} (via the exact same Morse theoretic argument), Kirby diagrams for non-orientable 4-manifolds are unique up to standard moves. 
\begin{theorem}\label{thm:kirbyunique}
Suppose $\mathcal{K}_0$ and $\mathcal{K}_1$ are Kirby diagrams for the 4-manifold $X^4$. Then $\mathcal{K}_0$ and $\mathcal{K}_1$ are related by a finite sequence of handle-slides and births/deaths of cancelling handle pairs.
\end{theorem}

\begin{example}
As an illustration of the techniques of this section, we will draw Kirby diagrams for $S^2\ttimes \mathbb{RP}^2$ and $\mathbb{RP}^4\#_{S^1}\mathbb{RP}^4$. We begin with the well known diagram for $D^2\ttimes \mathbb{RP}^2$ in Figure \ref{fig:fiberstructure}{c} (see \cite[Section 1.5]{akbulutbook}: it can be built with a single non-orientable $1$-handle, and a $2$-handle attached along a curve (in green) that runs across the $1$-handle twice, attached with framing as described below. We will adopt the convention that the feet of any non-orientable 1-handles in a diagram are identified via the identity map, which differs slightly from Akbulut's convention \cite{akbulutbook}. We handle the framings on $2$-handle curves in the same way as Akbulut, by indicating a framing relative to blackboard framing on each arc of a $2$-handle attaching circle minus the $1$-handles. Note that if a twist is pushed through a non-orientable $1$-handle, it becomes a twist of opposite sign.

The boundary of this disk bundle over $\mathbb{RP}^2$ is $S^2\ttimes S^1$, whose fibration structure is illustrated in Figure \ref{fig:fiberstructure}{a}. The exterior of the $2$-handle attaching curve is fibered by pairs of disks $D_\theta$ (one such disk is shaded blue in Figure \ref{fig:fiberstructure}{a}); the disks $D_\theta$ and $D_{\theta+\pi}$ glue together to form an annulus which meets the boundary of a neighbourhood of the 2-handle curve in two longitudes. Since we attach the $2$-handle with precisely the correct framing, we see the (twisted) $S^2$-bundle structure in $\partial (D^2\ttimes \mathbb{RP}^2)$.   

\begin{figure}[ht]
    \centering
     \scalebox{1}{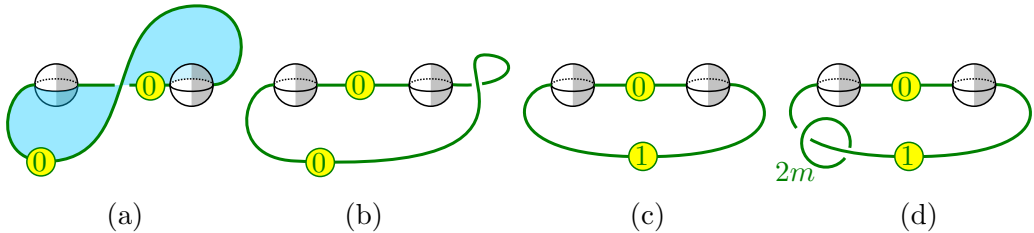}
     \put(-350,-15){(a)}
     \put(-260,-15){(b)}
         \put(-152,-15){(c)}
    \put(-50,-15){(d)}
    \caption{{\bf{(a):}} A disk making up a hemisphere of one sphere leaf of the fibration of $\partial( D^2\ttimes \mathbb{RP}^2)\cong S^2\ttimes S^1$. {\bf{(b)}}, {\bf{(c):}} A relative Kirby diagram for $D^2\ttimes \mathbb{RP}^2$. {\bf{(d):}} If $m$ is even, this is a Kirby diagram for $S^2\ttimes \mathbb{RP}^2$. If $m$ is odd, this is a Kirby diagram for $\mathbb{RP}^4\#_{S^1}\mathbb{RP}^4$.}%
    \label{fig:fiberstructure}
\end{figure}

A Kirby diagram of the double of this bundle may be obtained in the usual way, by simply adding a $0$-framed meridian to the $2$-handle attaching curve. To draw a diagram of $\mathbb{RP}^4\#_{S^1}\mathbb{RP}^4$, note that $\tau$ fixes $\gamma$ pointwise, but adds (or subtracts) $2$ to the framing. Thus we will add a $\pm2$-framed meridian, as illustrated in Figure \ref{fig:fiberstructure} (d). What remains in both cases is a $3$- and $4$-handle, which do not need to be specified. The map $\tau$ has order two, so the framing of the meridian only matters modulo 4. The reader may also wish to check this by handle slides.
\end{example}

\section{A non-orientable version of Waldhausen's Theorem}\label{sec:waldhausen}

In Section \ref{sec:LP}, we showed that trisection diagrams exist for non-orientable 4-manifolds. In this section, we characterize when a triple of curves on a non-orientable surface describes a trisection diagram; in the orientable setting, this characterization follows from Waldhausen's theorem.

\begin{theorem}[\cite{waldhausen}, see e.g.\ \cite{schleimer}]
Fix $g\ge k\ge 0$. Any two genus $g$ Heegaard splittings of $\#_k S^2\times S^1$ are equivalent.
\end{theorem}

Waldhausen actually first studied Heegaard splittings of $S^3$. The following well-known theorem of Haken, usually known as ``Haken's lemma" can then be used to extend the theorem to $\#_k S^2\times S^1$.

\begin{lemma}[\cite{haken}]\label{lem:haken}
Let $M$ be a (potentially non-orientable) 3-manifold containing an essential 2-sphere. Let $\Sigma$ be a Heegaard surface for $M$. Then there exists an essential 2-sphere $S$ in $M$ that intersects $\Sigma$ in a simple closed curve.
\end{lemma}

The following consequence of Waldausen's theorem is used frequently in the literature, especially in the classification of small-genus trisections \cite{jeffalexgenus2}.

\begin{figure}
    \centering
    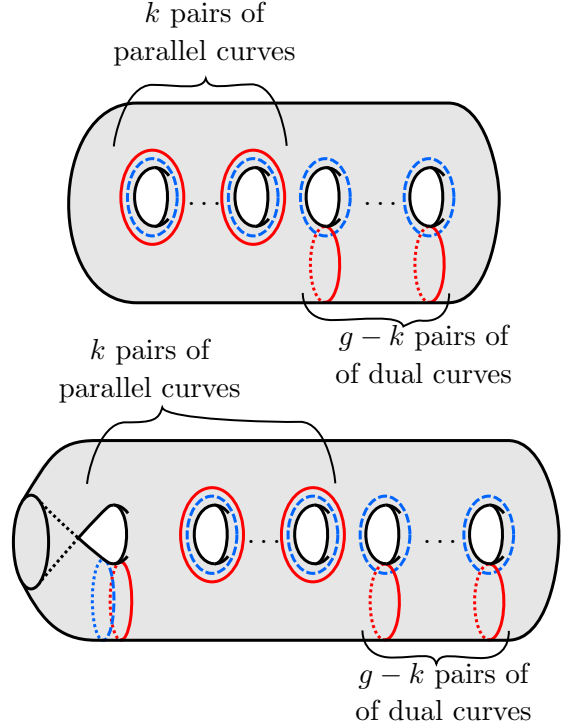
    \vspace{.1in}
    \caption{A genus $g$ Heegaard diagram of $\#_k S^2\times S^1$ (top) or $\#_k S^2\ttimes S^1$ (bottom). Here we draw one set of curves in red/solid and the other in blue/dashed. By Waldhausen's theorem and Theorem \ref{thm:waldhausen}, any genus $g$ Heegaard diagram $(\Sigma;\alpha,\beta)$ of $\#_kS^2\times S^1$ or $\#_kS^2\ttimes S^1$ is equivalent to one of these diagrams after an automorphism of $\Sigma$ and slides of $\alpha$ curves and slides of $\beta$ curves.}
    \label{fig:standardcurves}
\end{figure}

\begin{corollary}\label{cor:wald}
Let $\alpha,\beta,\gamma$ be cut systems of $g$ curves on an orientable genus $g$ surface $\Sigma$. Then $(\Sigma;\alpha,\beta,\gamma)$ is a trisection diagram if and only if each pair of $\alpha,\beta,\gamma$ is slide-equivalent to standard curves on $\Sigma$ (up to automorphism of $\Sigma$) as in the top of Figure \ref{fig:standardcurves}. We say that such a pair of curves is {\emph{standardizeable}}.
\end{corollary}

\begin{remark}\label{rem:simultaneous}
Corollary \ref{cor:wald} can be restated as, ``$(\Sigma;\alpha,\beta,\gamma)$ is a trisection diagram if and only if $(\alpha,\beta)$, $(\beta,\gamma)$, $(\gamma,\alpha)$ are each standardizeable." Note that although each pair of curves in a trisection diagram is standardizable, we do not expect different pairs to be {\emph{simultaneously}} standardizeable.
\end{remark}

Now we prove the analogous version of Waldhausen's theorem in the setting of non-orientable manifolds.

\begin{theorem}\label{thm:waldhausen}
Fix $g\ge k\ge 0$. Any genus $g$ Heegaard surface of $\#_k S^2\ttimes S^1$ is equivalent to the result of stabilizing the standard Heegaard surface $g-k$ times.
\end{theorem}

\begin{remark}
We remind the reader that $\#_k S^2\ttimes S^1$ admits a genus $k$ Heegaard splitting obtained by identifying two non-orientable genus $k$ handlebodies and gluing their boundaries by the identity map. The Heegaard surface of this splitting is called the {\emph{standard surface}} $\#_k S^2\ttimes S^1$; it is well-defined up to equivalence.
\end{remark}

Again, this theorem will follow from Waldhausen's work in $S^3$ via Haken's lemma.

\begin{proof}
If $k=0$, then the claim holds by Waldhausen's theorem. 

Now fix some $K\ge 1$ and suppose the claim holds whenever $k< K$, for all $g\ge k$. Write $\Sigma_{g,k}$ to indicate the result of stabilizing the standard surface in $\#_k S^2\ttimes S^1$ $g-k$ times.

Let $\Sigma$ be a genus $g$ Heegaard surface in $M:=\#_K S^2\ttimes S^1$. By Haken's lemma, there is an essential 2-sphere $S$ in $M$ that intersects $\Sigma$ in a simple closed curve. Let $(M',\Sigma')$ be the 3-manifold and Heegaard surface obtained by compressing $(M,\Sigma)$ along $S$.

\vspace{2mm}
\noindent{\bf{Case 1: $S$ is separating.}}
\vspace{2mm}

If $\Sigma$ is separating, then $M'$ is a disjoint union $M'=M_1\sqcup M_2$ with $M_i\cong \#_{k_i}S^2\times S^1$ or $M_i\cong\#_{k_i}S^2\ttimes S^1$, with $k_1+k_2=K$. Since $S$ is essential, $k_1,k_2\neq 0$. By Waldhausen's theorem or the inductive hypothesis, the corresponding components of $\Sigma'$ are equivalent to $\Sigma_{g_1,k_1}$ and $\Sigma_{g_2,k_2}$ for some $g_1,g_2$ with $g_1+g_2=G$. Then $(M,\Sigma)\cong(\#_KS^2\ttimes S^1,\Sigma_{G,K})$.

\vspace{2mm}
\noindent{\bf{Case 2: $S$ is non-separating.}}
\vspace{2mm}

If $\Sigma$ is nonseparating, then either $M'\cong \#_{K-1}S^2\times S^2$ or $M'\cong\#_{K-1}S^2\ttimes S^1$. By Waldhausen's theorem or the inductive hypothesis, $\Sigma'$ is equivalent to $\Sigma_{G-1,K-1}$. Then $(M,\Sigma)\cong(\#_KS^2\ttimes S^1,\Sigma_{G,K})$.

We thus conclude the claim holds for $k=K$, so Theorem \ref{thm:waldhausen} holds by induction.

\end{proof}

\begin{corollary}
Let $(\Sigma;\alpha,\beta)$ be a Heegaard diagram of a 3-manifold $M$. Then $M\cong\#_k S^2\ttimes S^1$ if and only if the following are satisfied:
\begin{itemize}
    \item $\Sigma$ is a non-orientable surface of genus $g\ge k$,
    \item $\alpha,\beta$ are slide-equivalent to the two sets of $g$ curves on $\Sigma$ indicated in the bottom of Figure \ref{fig:standardcurves} (up to automorphism of $\Sigma$).
\end{itemize}
\end{corollary}

\begin{corollary}
Let $\alpha,\beta,\gamma$ be cut systems of $g$ curves on a non-orientable genus $g$ surface $\Sigma$. Then $(\Sigma;\alpha,\beta,\gamma)$ is a trisection diagram if and only if each pair of $\alpha,\beta,\gamma$ is slide-equivalent to the standard curves in Figure \ref{fig:standardcurves} for some $k$ (up to automorphism of $\Sigma$). We say that such a pair of curves is {\emph{standardizeable}}; so $(\Sigma;\alpha,\beta,\gamma)$ is a trisection diagram if and only if $(\alpha,\beta)$, $(\beta,\gamma)$, $(\gamma,\alpha)$ are each standardizeable.
\end{corollary}

\begin{remark}
As in the orientable case (see Remark \ref{rem:simultaneous}), we do not expect the pairs $(\alpha,\beta)$, $(\beta,\gamma)$, $(\gamma,\alpha)$ to be {\emph{simultaneously}} standardizable even if $(\Sigma;\alpha,\beta,\gamma)$ is a non-orientable trisection diagram.
\end{remark}

Finally, we conclude the standard bijection between trisection diagrams and trisections.
\begin{theorem}
Let $X^4=X_1\cup X_2\cup X_3$ be a trisected, non-orientable 4-manifold. Then $\mathcal{T}=(X_1,X_2,X_3)$ determines a trisection diagram $(\Sigma;\alpha,\beta,\gamma)$ describing $\mathcal{T}$ that is well-defined up to automorphism of $\Sigma$ and slides of $\alpha,\beta,\gamma$. That is, there is a natural bijection \[\frac{\text{$\{$trisection diagrams$\}$}}{\text{surface automorphism, slides}}\leftrightarrow\frac{\text{$\{$trisected 4-manifolds$\}$}}{\text{diffeomorphism}}.\]
\end{theorem}
\section{Some example trisections}\label{sec:nonorientabletrisections}


In this section, we give some examples of non-orientable trisections and their diagrams. 

\begin{example}
Two simple examples of trisections of non-orientable 4-manifolds are described by the diagrams below. (The left and right boundary circles are identified by reflection through a horizontal axis, yielding closed non-orientable surfaces.) The only non-orientable 4-manifold admitting a genus one trisection is $S^3\ttimes S^1$. Indeed, there is only one non-separating curve with an annular neighbourhood on the Klein bottle, and so the only possible genus one trisection is given in Figure \ref{fig:basicexamples}{a}; one can check directly that this is $S^3\ttimes S^1$. The trisection genus of $\mathbb{RP}^4$ is equal to two, and a diagram is given in Figure \ref{fig:basicexamples}{b}. In Figure \ref{fig:getrp4}, we show that (b) is truly a diagram of $\mathbb{RP}^4$.

\begin{figure}[ht]
    \centering
    \includegraphics[width=0.7\textwidth]{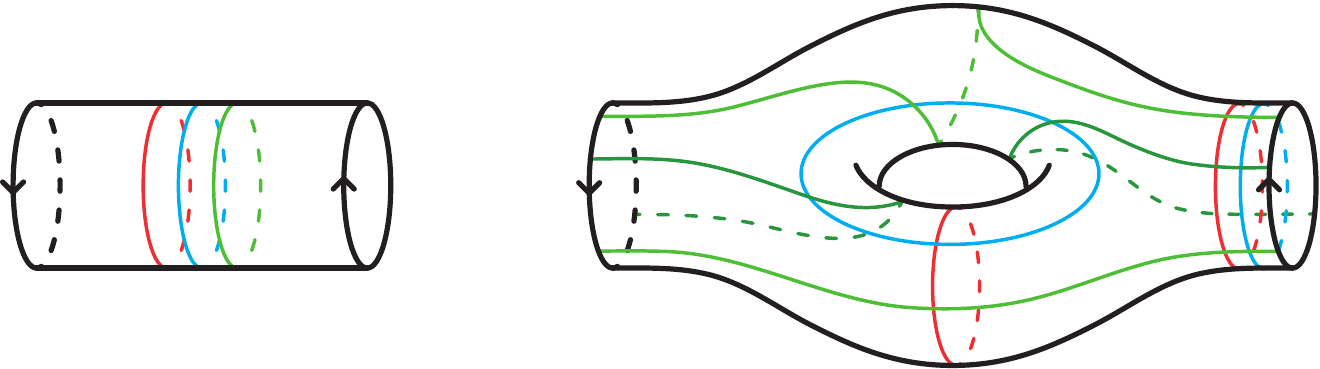}
    \put(-244,-15){(a)}
    \put(-80,-15){(b)}
    \caption{\textbf{(a):} a trisection diagram of $S^3\ttimes S^1$. \textbf{(b):} a trisection diagram for $\mathbb{RP}^4$.}%
    \label{fig:basicexamples}
\end{figure}

\begin{figure}
    \centering
    \scalebox{.7}{
\begingroup%
  \makeatletter%
  \providecommand\color[2][]{%
    \errmessage{(Inkscape) Color is used for the text in Inkscape, but the package 'color.sty' is not loaded}%
    \renewcommand\color[2][]{}%
  }%
  \providecommand\transparent[1]{%
    \errmessage{(Inkscape) Transparency is used (non-zero) for the text in Inkscape, but the package 'transparent.sty' is not loaded}%
    \renewcommand\transparent[1]{}%
  }%
  \providecommand\rotatebox[2]{#2}%
  \newcommand*\fsize{\dimexpr\f@size pt\relax}%
  \newcommand*\lineheight[1]{\fontsize{\fsize}{#1\fsize}\selectfont}%
  \ifx\svgwidth\undefined%
    \setlength{\unitlength}{463.36032657bp}%
    \ifx\svgscale\undefined%
      \relax%
    \else%
      \setlength{\unitlength}{\unitlength * \real{\svgscale}}%
    \fi%
  \else%
    \setlength{\unitlength}{\svgwidth}%
  \fi%
  \global\let\svgwidth\undefined%
  \global\let\svgscale\undefined%
  \makeatother%
  \begin{picture}(1,0.69118061)%
    \lineheight{1}%
    \setlength\tabcolsep{0pt}%
    \put(0,0){\includegraphics[width=\unitlength,page=1]{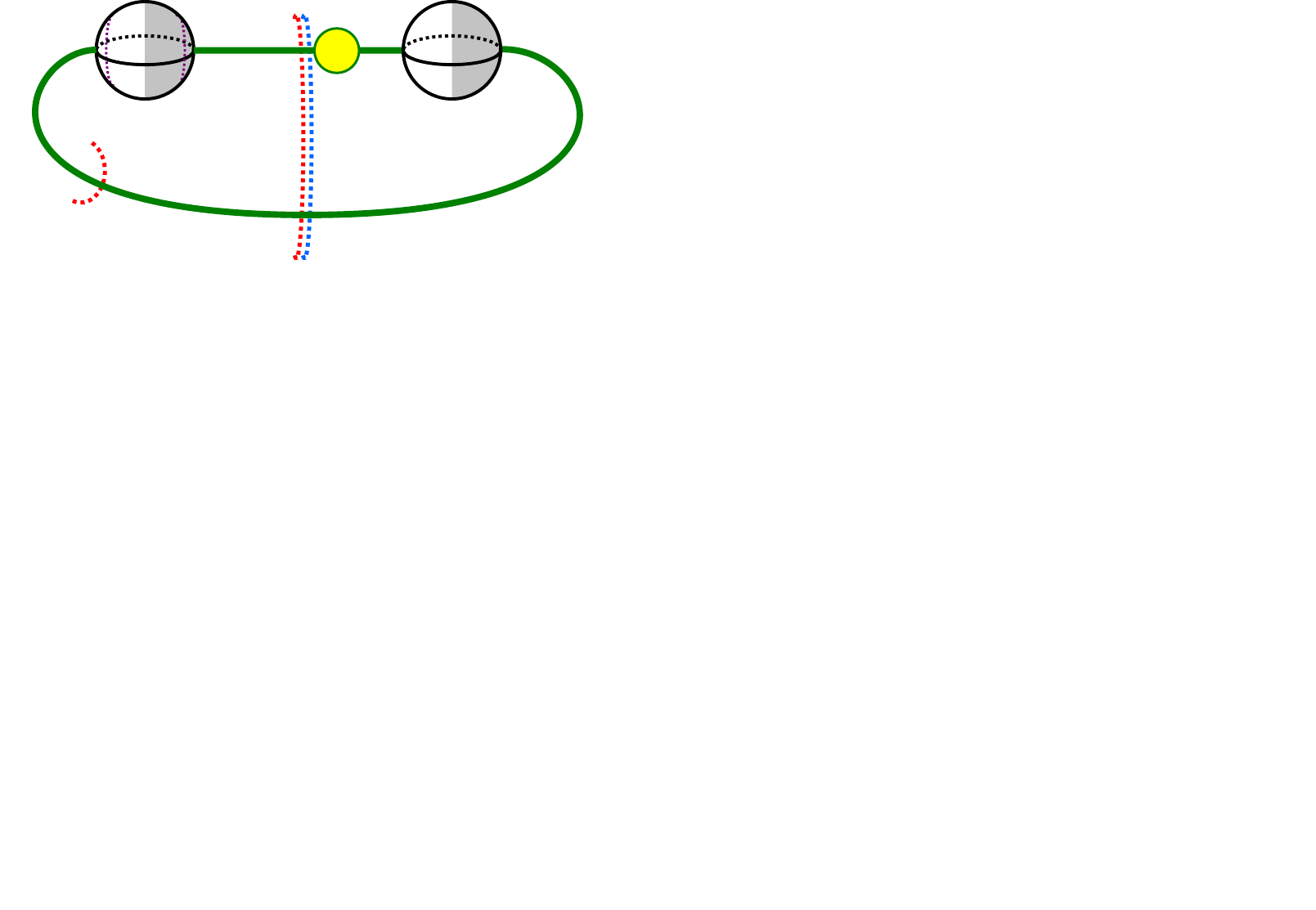}}%
    \put(0.24875802,0.64403687){\color[rgb]{0,0.50196078,0}\makebox(0,0)[lt]{\lineheight{1.25}\smash{\begin{tabular}[t]{l}$0$\end{tabular}}}}%
    \put(0,0){\includegraphics[width=\unitlength,page=2]{getrp4.pdf}}%
    \put(0.19371121,0.5195474){\color[rgb]{0,0.50196078,0}\makebox(0,0)[lt]{\lineheight{1.25}\smash{\begin{tabular}[t]{l}$1$\end{tabular}}}}%
    \put(0,0){\includegraphics[width=\unitlength,page=3]{getrp4.pdf}}%
  \end{picture}%
\endgroup%
}
    \caption{We start with the Kirby diagram of $\mathbb{RP}^4$ from Figure \ref{fig:fiberstructure}. We view this as a framed knot $K$ (the attaching circle of the 2-handle) inside $S^2\ttimes S^1$. We obtain a trisection diagram of $\mathbb{RP}^4$ as usual. {\bf{Top left:}} We find a Heegaard splitting of $\natural S^2\ttimes S^1$ in which $K$ lies in one handlebody as as a boundary-parallel curve.  {\bf{Top right:}} we project $K$ to the Heegaard surface with the correct framing to obtain one $\gamma$ curve. The other $\gamma$ curve is homologous to a combination of $\beta$ curves. {\bf{Bottom:}} we redraw our trisection diagram to see that we obtain the diagram of Figure \ref{fig:basicexamples}(b).}
    \label{fig:getrp4}
\end{figure}

\end{example}

\begin{example}\label{ex:cover}
By lifting each sector, a trisection of a closed non-orientable 4-manifold $X$ naturally lifts to a trisection of its orientation double cover, $\widetilde{X}$. Given a $(g;k)$-trisection for $X$, one obtains a $(2g-1;2k-1)$-trisection for $\widetilde{X}$. Moreover, the central surface $\widetilde{\Sigma}$ for $\widetilde{X}$ double covers the central surface $\Sigma$ for $X$, and so given a trisection diagram for $X$, we may easily draw the corresponding diagram for $\widetilde{X}$. Each curve in $\Sigma$ has an annular neighbourhood, and so will lift to two separate curves in $\widetilde{\Sigma}$; one curve among the lifts of each of the $\alpha$, $\beta$, and $\gamma$ cut systems will be homologically dependent, and so can be discarded. The $(3;1)$-trisection diagram of $S^4$ in Figure \ref{fig:doublecoverRP4} is the double cover of the standard trisection of $\mathbb{RP}^4$ in Figure \ref{fig:basicexamples}{b}. One can easily check that this trisection is handle slide diffeomorphic to the balanced stabilization of $S^4$ from \cite{gaykirby}. 

\vspace{3mm}
\begin{figure}[ht]
    \centering
    \includegraphics[width=0.7\textwidth]{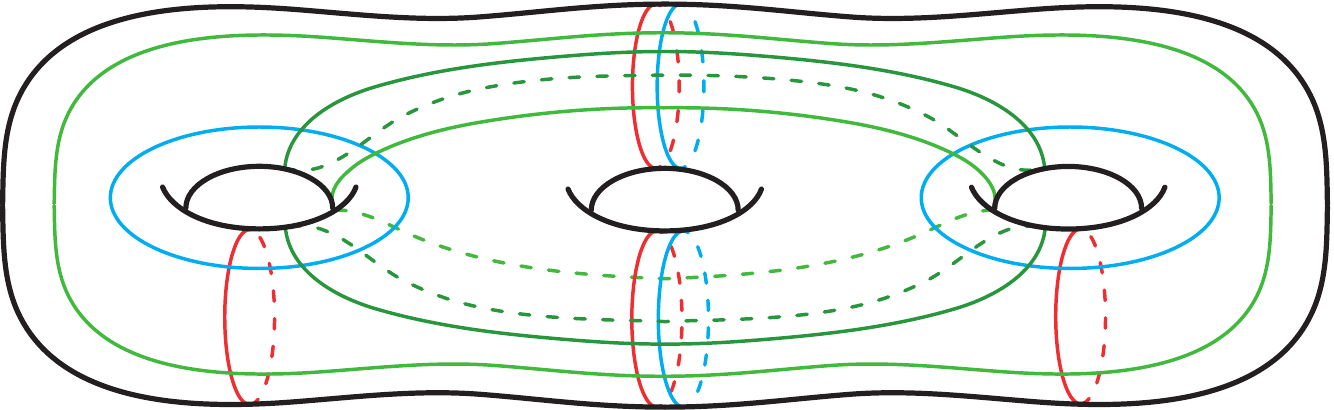}
    \caption{The (standard) $(3;1)$-trisection of $S^4$ induced by the trisection of $\mathbb{RP}^4$ in Figure \ref{fig:basicexamples}{b}. Each of the $\alpha$, $\beta$, and $\gamma$ cut systems contains a homologically dependent curve, which may be discarded. }%
    \label{fig:doublecoverRP4}
\end{figure}

\end{example}

\begin{remark}
A more interesting version of Example \ref{ex:cover} is obtained by considering Cappell-Shaneson homotopy 4-spheres \cite{cappell}, which are double covers of exotic $\mathbb{RP}^4$s. By lifting a trisection of such an exotic $\mathbb{RP}^4$, we obtain a trisection of a Cappell-Shaneson homotopy 4-sphere.

This is primarily interesting as a possible technique for standardizing these manifolds via extra symmetry of the trisection structure.

Secondarily, whether every trisection of $S^4$ is a stabilization of the $(0;0)$-trisection is an open problem related to the Andrews--Curtis Conjecture \cite{jeffalextrent}. A trisection of $S^4$ that double covers a trisection of an exotic $\mathbb{RP}^4$ at least cannot be equivariantly destabilized to the $(0;0)$-trisection.
\end{remark}

\begin{example}

For a slightly more interesting example, we consider the manifold $D^2\ttimes \mathbb{RP}^2\cup_{\tau^m}  D^2\ttimes \mathbb{RP}^2$ from Remark \ref{rem:diffeotopygroupS1xS2}, where $\tau$ is the twist map along a fiber in $S^2\ttimes S^1$. A Kirby diagram is illustrated in Figure \ref{fig:fiberstructure}, from which we can obtain the trisection diagram in Figure \ref{fig:twisteddouble}. 
These two trisections are minimal genus: the Euler characteristic of each of these manifolds is equal to 2, and so they cannot admit a trisection of genus less than 3.

In particular, this construction illustrates how one can obtain a trisection diagram from a Kirby diagram or vice versa.  This procedure is very well understood in the orientable category (see e.g.\ \cite{gaykirby} or \cite{jeffalextrent}), where it works in exactly the same way, so we omit detailed discussion.
%
\begin{figure}
    \centering
    \scalebox{.7}{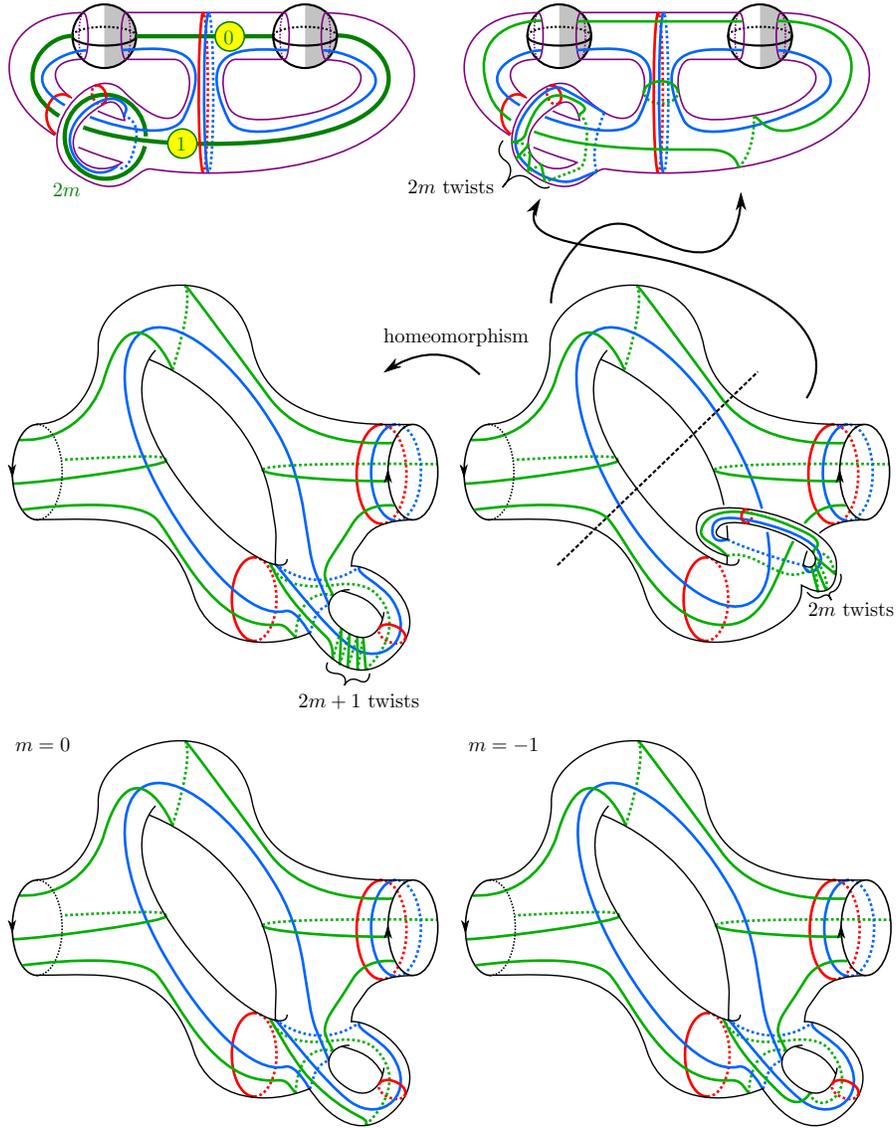}
    \caption{As in Figure \ref{fig:getrp4}, we start with a Kirby diagram consisting of a framed link $L$ in $S^2\ttimes S^1$. If $m$ is even, then this is a diagram of $S^2\ttimes \mathbb{RP}^2$; if $m$ is odd then this is a diagram of $\mathbb{RP}^4\#_{S^1}\mathbb{RP}^4$. {\bf{Top row:}} we obtain a trisection diagram for the described 4-manifold. {\bf{Middle row:}} We redraw the complete trisection diagram and simplify. {\bf{Bottom left:}} A trisection diagram of $S^2\ttimes \mathbb{RP}^2$. {\bf{Bottom right:}} A trisection diagram of $\mathbb{RP}^4\#_{S^1}\mathbb{RP}^4$.}
    \label{fig:twisteddouble}
\end{figure}

\end{example}

We end this section with a question: we will call a genus two trisection of a non-orientable 4-manifold \emph{standard} if it is either reducible or the trisection of $\mathbb{RP}^4$ above.

\begin{question}\label{genus2question}
Are all genus two non-orientable trisections standard? 
\end{question}

\begin{remark}
The unique non-orientable $(2;2)$-trisection is a reducible trisection of $\#_2 S^3\ttimes S^1$. Thus, Question \ref{genus2question} is specifically a question about (2;1)-trisections. Orientable (2;1)-trisections are all reducible \cite{jeffalexgenus2}; this follows from the fact that there are no nontrivial cosmetic surgeries in $S^2\times S^1$ \cite{gab86}. In contrast, $S^2\ttimes S^1$ admits a cosmetic surgery on a knot that wraps twice around the $S^1$ factor. This phenomenon yields the genus two non-orientable trisection of $\mathbb{RP}^4$ in Figure \ref{fig:basicexamples}(b). It is unlikely that Question \ref{genus2question} can be answered by lifting to the orientation double cover, as the lift of a $(2;1)$-non-orientable trisection is a $(3;1)$-orientable trisection, which are currently only conjecturally \cite{jeffspun} classified. Instead, one would likely need to understand all cosmetic surgeries on $S^2\ttimes S^1$ (which we suspect includes a unique nontrivial surgery).
\end{remark}

\section{Non-orientable 4-manifolds with boundary}\label{sec:relative}

In this section, we will discuss relative trisections of 4-manifolds with boundary. To understand these decompositions (even for orientable 4-manifolds), we must understand {\emph{compression bodies}}. We refer the reader to \cite{castro} or \cite{casgaypin18} for more detailed treatment of relative trisections of orientable 4-manifolds; we follow the notation of \cite{castro} here. In particular, \cite{casgaypin18} describes the case in which the 4-manifold has many boundary components, while \cite{castro} might be called, ``the springboard into the sea of relative trisections."\footnote{Personal communication with N. A. Castro, 29 Sep. 2020 (possibly a joke).}

\subsection{Definition}

\begin{figure}
    \centering
\begingroup%
  \makeatletter%
  \providecommand\color[2][]{%
    \errmessage{(Inkscape) Color is used for the text in Inkscape, but the package 'color.sty' is not loaded}%
    \renewcommand\color[2][]{}%
  }%
  \providecommand\transparent[1]{%
    \errmessage{(Inkscape) Transparency is used (non-zero) for the text in Inkscape, but the package 'transparent.sty' is not loaded}%
    \renewcommand\transparent[1]{}%
  }%
  \providecommand\rotatebox[2]{#2}%
  \newcommand*\fsize{\dimexpr\f@size pt\relax}%
  \newcommand*\lineheight[1]{\fontsize{\fsize}{#1\fsize}\selectfont}%
  \ifx\svgwidth\undefined%
    \setlength{\unitlength}{282.59015902bp}%
    \ifx\svgscale\undefined%
      \relax%
    \else%
      \setlength{\unitlength}{\unitlength * \real{\svgscale}}%
    \fi%
  \else%
    \setlength{\unitlength}{\svgwidth}%
  \fi%
  \global\let\svgwidth\undefined%
  \global\let\svgscale\undefined%
  \makeatother%
  \begin{picture}(1,0.6145742)%
    \lineheight{1}%
    \setlength\tabcolsep{0pt}%
    \put(0,0){\includegraphics[width=\unitlength,page=1]{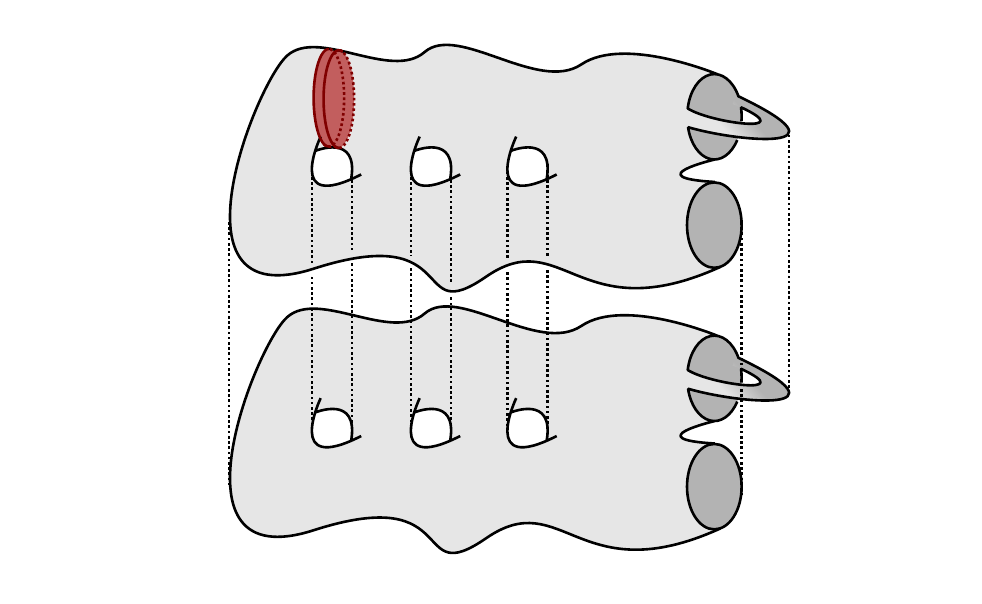}}%
    \put(-0.00252526,0.30425407){\color[rgb]{0,0,0}\makebox(0,0)[lt]{\lineheight{1.25}\smash{\begin{tabular}[t]{l}$\Sigma\times[0,1]$\end{tabular}}}}%
    \put(0.43705117,0.00717467){\color[rgb]{0,0,0}\makebox(0,0)[lt]{\lineheight{1.25}\smash{\begin{tabular}[t]{l}$\partial_-C$\end{tabular}}}}%
    \put(0.01548893,0.50777802){\color[rgb]{0,0,0}\makebox(0,0)[lt]{\lineheight{1.25}\smash{\begin{tabular}[t]{l}2-handles\end{tabular}}}}%
    \put(0,0){\includegraphics[width=\unitlength,page=2]{compressionbody.pdf}}%
    \put(0.43705117,0.5914607){\color[rgb]{0,0,0}\makebox(0,0)[lt]{\lineheight{1.25}\smash{\begin{tabular}[t]{l}$\partial_+C$\end{tabular}}}}%
    \put(0.81207439,0.32841475){\color[rgb]{0,0,0}\makebox(0,0)[lt]{\lineheight{1.25}\smash{\begin{tabular}[t]{l}$(\partial \Sigma)\times[0,1]$\end{tabular}}}}%
  \end{picture}%
\endgroup%

    \caption{A compression body $C$ on a surface $\Sigma$. The boundary of $C$ decomposes as $\partial_-C\cup(\boundary\Sigma)\times[0,1]\cup\partial_+ C$.}
    \label{fig:compressionbody}
\end{figure}

\begin{definition}
A {\emph{compression body}} is a 3-manifold $C$ of the form $C=(\Sigma\times[0,1]\cup(2$-handles glued along curves in $\Sigma\times\{1\})$, where $\Sigma$ is a connected surface with nonempty boundary. See Figure \ref{fig:compressionbody}. We may call $C$ a {\emph{compression body on $\Sigma$}}. Note that $C$ is orientable if and only if the surface $\Sigma$ is orientable. We have \[\boundary C=\boundary_-C\cup(\boundary \Sigma\times(0,1))\cup\boundary_+C,\]
where \[\boundary_-C:=\Sigma\times\{0\}\] and \[\boundary_+C:=\boundary C\setminus\left(\boundary_-C\cup(\boundary\Sigma\times(0,1))\right).\]
We call $\boundary_-C$ and $\boundary_+C$ the {\emph{negative boundary}} and {\emph{positive boundary}} of $C$, respectively. If $\Sigma$ is not orientable, then neither is $\boundary_- C$. However, $\boundary_+ C$ may be orientable regardless. 
\end{definition}

\begin{exercise}
Let $C$ be a compression body. Then 
\[C\times I\cong \begin{cases}\natural_k B^3\times S^1&\text{if $C$ is orientable}\\
\natural_k B^3\ttimes S^1&\text{if $C$ is non-orientable},\end{cases}\] 
where 
\[k=1-\frac{\chi(\boundary_+ C)+\chi(\boundary_-C)}{2}.\]
\end{exercise}

\begin{definition}
Let $C$ be a compression body on a surface $\Sigma$, and set $Z:=C\times I$. We decompose $\partial Z=\boundaryin Z\cup\boundaryout Z$ into two pieces, where 
\[\boundaryin Z:=(C\times \{0\})\cup(\partial_-C\times I)\cup(C\times\{1\})\] and 
\[\boundaryout Z:=(\partial \Sigma\times I\times I)\cup(\partial_+C\times I).\]
We call $\boundaryin Z$ and $\boundaryout Z$ the {\emph{interior boundary}} and {\emph{outer boundary}} of $Z$, respectively.  See Figure \ref{fig:boundaryinout}.

\begin{figure}
    \centering
\begingroup%
  \makeatletter%
  \providecommand\color[2][]{%
    \errmessage{(Inkscape) Color is used for the text in Inkscape, but the package 'color.sty' is not loaded}%
    \renewcommand\color[2][]{}%
  }%
  \providecommand\transparent[1]{%
    \errmessage{(Inkscape) Transparency is used (non-zero) for the text in Inkscape, but the package 'transparent.sty' is not loaded}%
    \renewcommand\transparent[1]{}%
  }%
  \providecommand\rotatebox[2]{#2}%
  \newcommand*\fsize{\dimexpr\f@size pt\relax}%
  \newcommand*\lineheight[1]{\fontsize{\fsize}{#1\fsize}\selectfont}%
  \ifx\svgwidth\undefined%
    \setlength{\unitlength}{283.58031265bp}%
    \ifx\svgscale\undefined%
      \relax%
    \else%
      \setlength{\unitlength}{\unitlength * \real{\svgscale}}%
    \fi%
  \else%
    \setlength{\unitlength}{\svgwidth}%
  \fi%
  \global\let\svgwidth\undefined%
  \global\let\svgscale\undefined%
  \makeatother%
  \begin{picture}(1,0.65770515)%
    \lineheight{1}%
    \setlength\tabcolsep{0pt}%
    \put(0,0){\includegraphics[width=\unitlength,page=1]{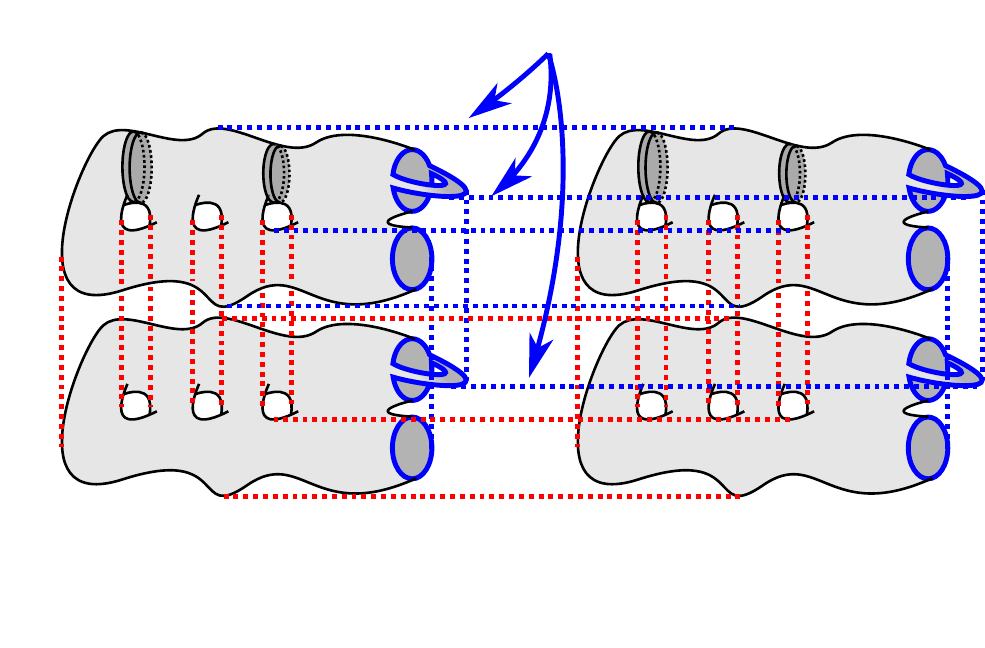}}%
    \put(0.286417,0.62609961){\color[rgb]{0,0,1}\makebox(0,0)[lt]{\lineheight{1.25}\smash{\begin{tabular}[t]{l}$\boundaryout Z=(\partial\Sigma\times I\times I)\cup(\partial_+ C\times I)$\end{tabular}}}}%
    \put(0.21455269,0.00817127){\color[rgb]{1,0,0}\makebox(0,0)[lt]{\lineheight{1.25}\smash{\begin{tabular}[t]{l}$\boundaryin Z=(C\times\{0\})\cup(\partial_- C\times I)\cup(C\times\{1\})$\end{tabular}}}}%
    \put(0,0){\includegraphics[width=\unitlength,page=2]{boundaryinout.pdf}}%
  \end{picture}%
\endgroup%

    \caption{A schematic of the decomposition $\boundary Z=\boundaryin Z\cup\boundaryout Z$.}
    \label{fig:boundaryinout}
\end{figure}

The interior boundary of $Z$ comes with a natural (generalized) Heegaard splitting $\boundaryin Z=Y_0^-\cup Y_0^+$, where \[Y_0^-=(C\times \{0\})\cup(\partial_-C\times(0,1/2]), \hspace{.5in}Y_0^+=(\partial_-C\times[1/2,1))\cup(C\times\{1\}).\] The Heegaard surface of this splitting is $\partial_-C\times\{1/2\}$. Any Heegaard splitting of $\boundaryin Z$ obtained by stabilizing $(Y_0^-,Y_0^+)$ some (possibly zero) number of times is called a {\emph{standard splitting}} of $\boundaryin Z$. Note that standard splitting $(Y^-,Y^+)$ of $\boundaryin Z$ is determined up to isotopy by the genus of $Y^-\cap Y^+$.
\end{definition}

\begin{notation}
For the rest of this section, $C$ will always refer to a compression body on a surface $\Sigma$ and $Z$ will always refer to $C\times I$. Whenever we write $\boundaryin Z=Y^-\cup Y^+$, we implicitly mean that $(Y^-,Y^+)$ is a standard splitting of $\boundaryin Z$.
\end{notation}

\begin{figure}
    \centering
    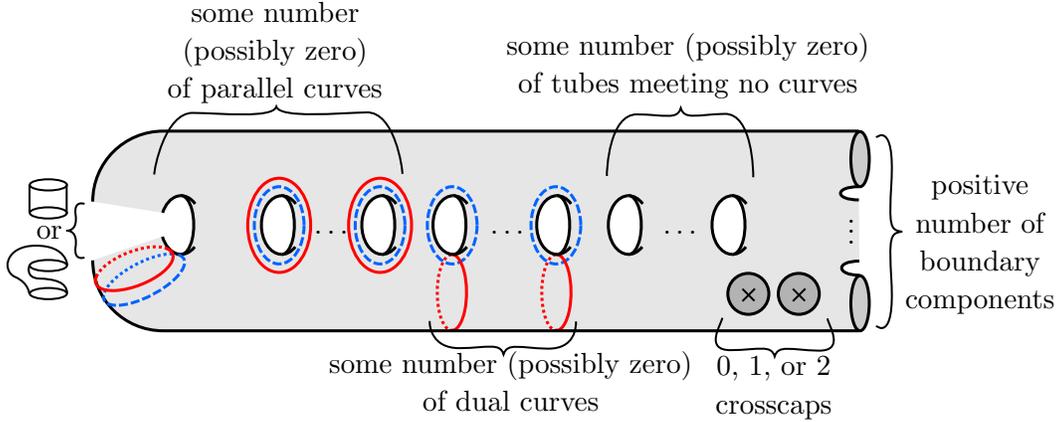
    \caption{A standard pair (of sets of curves) on a possibly non-orientable surface $\Sigma$ with nonempty boundary. One set of curves is red/solid, the other is blue/dashed.}
    \label{fig:reldiagramdef}
\end{figure}

\begin{definition}
Let $\Sigma$ be a (possibly non-orientable) surface of genus $g\ge 0$ and $b>0$ boundary components, and let $\alpha,\beta$ be sets of disjoint simple closed curves on $\Sigma$. We say that $(\alpha,\beta)$ is a {\emph{standard}} pair if 
there exists a homeomorphism of $\Sigma$ taking $\alpha$ to the red (solid) curves and $\beta$ to the blue (dashed) curves in Figure \ref{fig:reldiagramdef}. 

If there exist sets of curves $\alpha',\beta'$ in $\Sigma$ with $\alpha$ slide-equivalent to $\alpha'$, $\beta$ slide-equivalent to $\beta'$, and $(\alpha',\beta')$ a standard pair, then we say that $(\alpha,\beta)$ is a {\emph{slide-standard}} pair.
\end{definition}

\begin{proposition}\label{prop:standardsplit}
Let $(\alpha,\beta)$ be a slide-standard pair in a surface $\Sigma$. Then $(\Sigma;\alpha,\beta)$ determines a standard splitting $\boundaryin Z=Y^-\cup Y^+$. 
\end{proposition}

\begin{proof}
Let $V:=\Sigma\times I\cup H_\alpha\cup H_\beta$, where $H_\alpha$ are $2$-handles attached along $\alpha\times\{0\}$ and $H_\beta$ are $2$-handles attached along $\beta\times\{1\}$. Let $Y^-:=H_\alpha\cup(\Sigma\times[0,1/2])$ and $Y^+=(\Sigma\times[1/2,1])\cup H_\beta$.

Then $V\cong \boundaryin Z$ for $Z=C\times I$ and $C$ a compression body on a surface $\Sigma'$ with $\chi(\Sigma')=\chi(\Sigma)+2n,|\boundary\Sigma'|=|\boundary \Sigma|$, and $\Sigma'$ is orientable if and only if $\Sigma$ orientable, where $n$ is the number of pairs of dual $\alpha,\beta$ curves. The pair $(Y^-,Y^+)$ is obtained by stabilizing $(Y_0^-,Y_0^+)$ a total of $n$ times. Thus, $(Y^-,Y^+)$ is a standard splitting of $\boundaryin Z$.
\end{proof}

We are now prepared to define a relative trisection.

\begin{definition}\label{def:reltri}
Let $X^4$ be a compact 4-manifold with connected, nonempty boundary. A {\emph{relative trisection}} of $X$ is a decomposition $X=X_1\cup X_2\cup X_3$ such that:
\begin{itemize}
    \item There exist diffeomorphisms $\phi_i:X_i\to Z$ with $\phi_i(X_i\cap\boundary X)=\boundaryout Z$,
    \item For each $i$, $\phi_i(X_i\cap X_{i-1})=Y^-$ and $\phi_i(X_{i}\cap X_{i+1})=Y^+$.
\end{itemize}

This definition can be extended to include manifolds with disconnected boundary \cite{casgaypin18}, but for simplicity we ignore this case.
\end{definition}

\subsection{Gluing}

In order to glue relatively trisected 4-manifolds, we must understand the boundary of a trisected 4-manifold.

\begin{proposition}
A relative trisection $X=X_1\cup X_2\cup X_3$ induces an open book decomposition on $\boundary X$.
\end{proposition}

\begin{proof}
Observe that $L:=\boundary(X_1\cap X_2\cap X_3)$ is a link in $\boundary X$. From Definition \ref{def:reltri}, we have $\phi_i(X_i\cap\boundary X)=\boundaryout Z=(\partial\Sigma\times I\times I)\cup(\boundary_+C\times I)$. Then $X_i\cap\boundary X$ is a thickened Seifert surface for $L$ along with a tubular neighborhood of $L$. We conclude that for some tubular neighborhood $\nu(L)$, $\partial X\setminus\nu(L)$ is the total space of a bundle of surfaces over $S^1$, where $X_j\cap X_{j+1}\cap(\partial X\setminus\nu(L))$ is a fiber for each $j$. That is, $L$ is the binding of an open book on $X$ which includes the pages $X_j\cap X_{j+1}\cap(\partial X)$ for $j=1,2,3$. 
\end{proof}

If $X$ is orientable, then given an open book $\mathcal{O}$ on $\partial X$, there exists a relative trisection of $X$ that induces $\mathcal{O}$ \cite{gaykirby,castro}.

\begin{theorem}[\cite{gaykirby}]\label{reltriunique}
Let $\mathcal{T}_1$ and $\mathcal{T}_2$ be relative trisections of a 4-manifold $X^4$ inducing isotopic open books on $\partial X^4$. Then $\mathcal{T}_1$ and $\mathcal{T}_2$ agree after a finite sequence of interior stabilizations in each.
\end{theorem}

Theorem \ref{reltriunique} is stated only in the orientable setting, but the proof carries through for non-orientable manifolds. Alternatively, one can prove Theorem \ref{reltriunique} (or similarly for closed trisections) via (a relative version of) Theorem \ref{thm:kirbyunique} by first converting each of $\mathcal{T}_1,\mathcal{T}_2$ to a relative handle decomposition.

\begin{remark}
Any two relative trisections of an {\emph{orientable}} 4-manifold $X$ are related by interior stabilization \cite{gaykirby}, relative double twists \cite{CIMT}, and relative stabilization \cite{castro}, in that order. (See the cited papers for definitions of these moves, if desired). This relies heavily on classification of open books on orientable 3-manifolds, especially work of Giroux--Goodman. To prove an analogous theorem for non-orientable relative trisections, one would need to find a complete set of moves relating any two open books of a non-orientable 3-manifold. As observed by Ozbagci \cite{ozbagci}, two open books of a non-orientable 3-manifold are not necessarily related by Hopf stabilization -- so we conclude that two relative trisections of a non-orientable 3-manifold are not necessarily related by a sequence of interior and relative stabilization.

\end{remark}

Perhaps the most important theorem about relative trisections is the following gluing theorem, proved by Castro in his thesis \cite{castro}.

\begin{theorem}[\cite{castro}]
Let $\T$ and $\T'$ be relative trisections of 4-manifolds $X$ and $X'$, respectively. Let $\mathcal{O}$ and $\mathcal{O}'$ denote the open books of $\boundary X$ and $\boundary X'$ (respectively) induced by $\T$ and $\T'$. Suppose there exists a diffeomorphism $f:\boundary X\to\boundary X'$ and that moreover $f(\mathcal{O})$ is isotopic to $\mathcal{O}'$. Then there is a naturally induced trisection $\T\cup\T'$ of $X\cup_f X'$. 
\end{theorem}
Of course, Castro stated this theorem only in the case that $X$ and $X'$ are orientable, but his proof holds in the non-orientable setting verbatim.

\subsection{Diagrams}

Relative trisections are also described by relative trisection diagrams -- see \cite{reldiagrams} in the orientable setting.
\begin{definition}
A {\emph{relative trisection diagram}} is a tuple $(\Sigma;\alpha,\beta,\gamma)$ where $\Sigma$ is a (possibly non-orientable) connected surface with nonempty boundary, $\alpha,\beta,\gamma$ are sets of disjoint simple closed curves in $\Sigma$, and each of $(\alpha,\beta),(\beta,\gamma),(\gamma,\alpha)$ is a slide-standard pair. 
\end{definition}

\begin{proposition}\label{prop:reldiagrams}
A relative trisection diagram $\mathcal{D}=(\Sigma;\alpha,\beta,\gamma)$ determines a relatively trisected 4-manifold $X=X_1\cup X_2\cup X_3$ up to diffeomorphism.
\end{proposition}

\begin{figure}
    \centering
     \scalebox{.9}{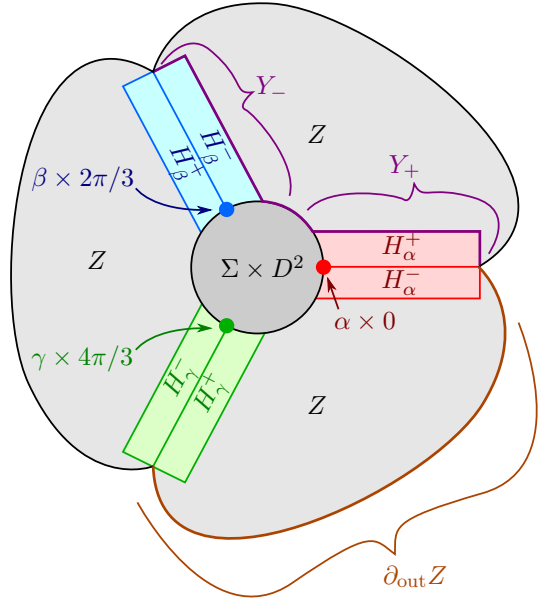}
    \caption{The recipe for building a 4-manifold with boundary from a relative trisection diagram, as described in Proposition \ref{prop:reldiagrams}.}
    \label{fig:buildrelmfd}
\end{figure}

Castro--Gay--Pinz\'on-Caicedo proved Proposition \ref{prop:reldiagrams} (and the ensuing Proposition \ref{prop:getreldiagram} and Corollary \ref{cor:reldiagrams}) for orientable trisection diagrams \cite{reldiagrams}. The proof follows in exactly the same manner for non-orientable diagrams, but we repeat the construction here. Alternatively, a reader familiar with trisections may refer to Figure \ref{fig:buildrelmfd}.

\begin{proof}
The diagram $\mathcal{D}$ immediately determines the {\emph{spine}} of a trisection: let $W$ be the 4-manifold with boundary obtained from $\Sigma\times D^2$ by attaching $2$-handles to $\Sigma\times S^1$ along $\alpha\times\{0\},\beta\times\{2\pi/3\},\gamma\times\{4\pi/3\}$. Specifically, let $\alpha_i$ be a curve in $\alpha$. The annular neighborhood of $\alpha_i\times\{0\}$ in $\Sigma\times\{0\}$ induces a framing on $\alpha_i\times\{0\}$ in $\boundary W$. Attach a 2-handle to $\alpha_i\times\{0\}$ via this framing. Repeat similarly for each curve in $\alpha,\beta,\gamma$. Let $H_\alpha,H_\beta,H_\gamma$ respectively denote the 2-handles attached along $\alpha\times\{0\},\beta\times\{2\pi/3\},\gamma\times\{4\pi/3\}$.

A 2-handle $H$ in $H_\alpha,H_\beta,H_\gamma$ is abstractly a $D^2\times D^2$ attached along $(\boundary D^2)\times D^2$. Write $D^2$ as the union of two half-disks $D^-,D^+$ along an equator, so the 2-handle $H$ can be expressed as $H^-\cup H^+$, where $H^-=(D^2\times D^-)$ and $H^+=(D^2\times D^+)$. If the attaching circle $\eta$ of $H$ is in $\Sigma\times\{\theta\}$, take $D^-\cap D^+$ to also lie in $\Sigma\times\{\theta\}$ with $D^-=(\nu(\eta)\cap(\Sigma\times\{\theta'\mid \theta'\le\theta\})$ and $D^+=(\nu(\eta)\cap(\Sigma\times\{\theta'\mid \theta'\le\theta\})$. Write $H^-_*$ (resp. $H^+_*$) to denote the $H^-$ (resp. $H^+$) part of each handle in $H_*$.

Now consider the following submanifold $V$ of $\boundary W$:
\[V:=\left((\Sigma\times[0,2\pi/3])\cup H^+_\alpha\cup H^-_\beta\right)\cap\boundary W. \]
    As in Proposition \ref{prop:standardsplit}, $V\cong \boundaryin Z$, in which $\alpha$ and $\beta$ have determined a standard splitting $(Y^-,Y^+)$ of determined genus. Therefore, there is a unique way to glue a thickened compression body $Z_1$ to $W$ along $V$.
    
Similarly, we glue thickened compression bodies $Z_2,Z_3$ (respectively) to
\[\left((\Sigma\times[2\pi/3,4\pi/3])\cup H^+_\beta\cup H^-_\gamma\right)\cap\boundary W. \] and \[\left((\Sigma\times[4\pi/3,2\pi]\cup H^+_\gamma\cup H^-_\alpha\right)\cap\boundary W. \]

The result is our desired relatively trisected 4-manifold. To see the relative trisection, set
\begin{align*}
    X_1&=\Sigma\times(\{0\le\theta\le2\pi/3\}\subset D^2)\cup Z_1,\\
     X_2&=\Sigma\times(\{2\pi/3\le\theta\le4\pi/3\}\subset D^2)\cup Z_2,\\
      X_3&=\Sigma\times(\{4\pi/3\le\theta\le2\pi\}\subset D^2)\cup Z_3.
\end{align*}

\end{proof}

The converse of Proposition \ref{prop:reldiagrams} is simpler.
\begin{proposition}\label{prop:getreldiagram}
Let $X=X_1\cup X_2\cup X_3$ be a relatively trisected 4-manifold. Then the relative trisection of $X$ determines a relative trisection diagram $(\Sigma;\alpha,\beta,\gamma)$ up to slides of $\alpha,\beta,\gamma$ and automorphism of $\Sigma$.
\end{proposition}
\begin{proof}
Let $\Sigma=X_1\cap X_2\cap X_3$. Recall there is a diffeomorphism $\phi_i:X_i\to Z$ with $\phi_1(X_i\cap X_{i+1})=Y^+$. Moreover, $Y^+$ can be obtained from $\phi_i(\Sigma\times I)$ by attaching 2-handles to $\Sigma\times\{1\}$. Let the attaching curves of these 2-handles, viewed as curves in $\Sigma$, be $\alpha,\beta,\gamma$ for $i=1,2,3$, respectively.

We have $\phi_1(X_3\cap X_{1})=Y^-$ and $\phi_1(X_1\cap X_{2})=Y^+$. Since $(Y^-,Y^+)$ is a standard splitting of $\boundaryin Z$, we find that $(\gamma,\alpha)$ are slide-standard. By similarly considering $\phi_2,\phi_3$, we find that $(\alpha,\beta)$ and $(\beta,\gamma)$ are slide standard. Therefore, $(\Sigma;\alpha,\beta)$ is a relative trisection diagram. 
\end{proof}

\begin{remark}
Proposition \ref{prop:getreldiagram} illustrates why the definition of a relative trisection is so much longer than the definition of a trisection (and hence why Section \ref{sec:relative} is so much longer than previous sections). Pragmatically, we want a relative trisection to determine a diagram uniquely up to some kind of equivalence. The fact that closed trisections determine unique trisection diagrams (up to automorpishm and slides) follows from Waldhausen's theorem or Theorem \ref{thm:waldhausen} applied to the boundary of each sector of the trisection. However, the sectors of a relative trisection do not have induced Heegaard splittings on their boundaries: they have only a splitting of the portion of the sector boundary contained in the interior of the manifold. If we indicate this partial splitting via curves on a central surface, then we cannot use \cite{laupoe72} or Theorem \ref{thm:NOLP} to glue the sectors together, since we have not specified their whole boundaries. Moreover, even if we somehow knew how to build $X$ from this ``diagram," we could not hope that the diagram would be determined uniquely by the trisection up to automorphism/slides.

Thus, the definition of a relative trisection is actually quite subtle. One may be tempted to think that if a decomposition of $X^4=X_1\cup X_2\cup X_3$ into three handlebodies with $X_i\cap X_j$ a handlebody and $X_1\cap X_2\cap X_3$ a neatly embedded surface induces an open book on $\partial X^4$, then this decomposition is a relative trisection. In fact, this definition would be sufficient for the purpose of gluing two relative trisections and obtaining a closed trisection, but we would not be able to use $(X_1,X_2,X_3)$ to produce a surface diagram that uniquely determines $X$ (and the trisection $(X_1,X_2,X_3)$) up to diffeomorphism. 
\end{remark}

Finally, we remark that the algorithms of Propositions \ref{prop:reldiagrams} and \ref{prop:getreldiagram} are clearly inverses. We therefore conclude the following corollary.

\begin{corollary}\label{cor:reldiagrams}
There is a natural bijection

\[\frac{\text{$\{$relative trisection diagrams$\}$}}{\text{surface automorphism, slides}}\leftrightarrow\frac{\text{$\{$relatively trisected 4-manifolds$\}$}}{\text{diffeomorphism}}.\]
\vspace{2mm}
\end{corollary}

\subsection{Examples}

\begin{figure}
    \centering
\begingroup%
  \makeatletter%
  \providecommand\color[2][]{%
    \errmessage{(Inkscape) Color is used for the text in Inkscape, but the package 'color.sty' is not loaded}%
    \renewcommand\color[2][]{}%
  }%
  \providecommand\transparent[1]{%
    \errmessage{(Inkscape) Transparency is used (non-zero) for the text in Inkscape, but the package 'transparent.sty' is not loaded}%
    \renewcommand\transparent[1]{}%
  }%
  \providecommand\rotatebox[2]{#2}%
  \newcommand*\fsize{\dimexpr\f@size pt\relax}%
  \newcommand*\lineheight[1]{\fontsize{\fsize}{#1\fsize}\selectfont}%
  \ifx\svgwidth\undefined%
    \setlength{\unitlength}{381.7013718bp}%
    \ifx\svgscale\undefined%
      \relax%
    \else%
      \setlength{\unitlength}{\unitlength * \real{\svgscale}}%
    \fi%
  \else%
    \setlength{\unitlength}{\svgwidth}%
  \fi%
  \global\let\svgwidth\undefined%
  \global\let\svgscale\undefined%
  \makeatother%
  \begin{picture}(1,0.25298715)%
    \lineheight{1}%
    \setlength\tabcolsep{0pt}%
    \put(0,0){\includegraphics[width=\unitlength,page=1]{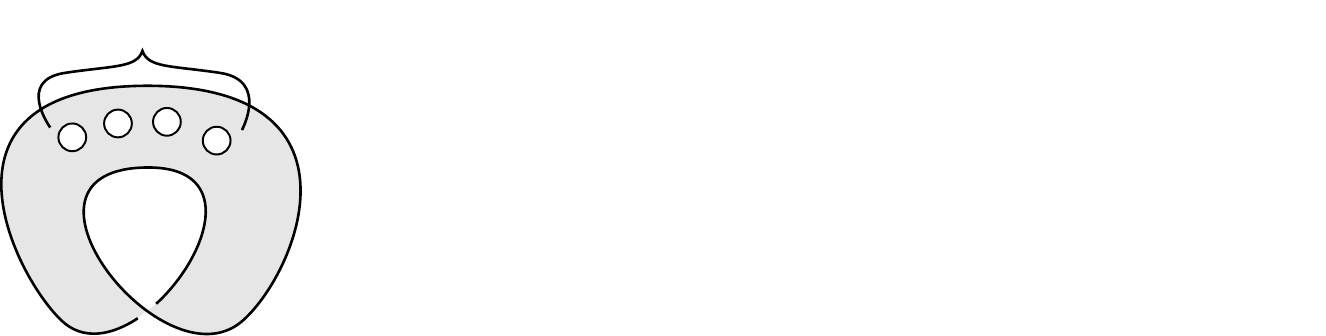}}%
    \put(0.06737176,0.22759503){\color[rgb]{0,0,0}\makebox(0,0)[lt]{\lineheight{1.25}\smash{\begin{tabular}[t]{l}$n\ge 0$\end{tabular}}}}%
    \put(0,0){\includegraphics[width=\unitlength,page=2]{smallrelative.pdf}}%
    \put(0.30455047,0.22759503){\color[rgb]{0,0,0}\makebox(0,0)[lt]{\lineheight{1.25}\smash{\begin{tabular}[t]{l}$n>0$\end{tabular}}}}%
    \put(0,0){\includegraphics[width=\unitlength,page=3]{smallrelative.pdf}}%
    \put(0.53843657,0.22759503){\color[rgb]{0,0,0}\makebox(0,0)[lt]{\lineheight{1.25}\smash{\begin{tabular}[t]{l}$n>0$\end{tabular}}}}%
    \put(0,0){\includegraphics[width=\unitlength,page=4]{smallrelative.pdf}}%
  \end{picture}%
\endgroup%

    \caption{Some small valid relative trisection diagrams. From left to right, these describe relative trisections of: $\natural_{n+1} B^3\ttimes S^1$, $\natural_{n+1}B^3\ttimes S^1$, $((S^3\ttimes S^1)\setminus\mathring{B}^4)\natural(\natural_{n-1}B^3\ttimes S^1)$, $D^2\ttimes \mathbb{RP}^2$.}
    \label{fig:smallrelative}
\end{figure}

In this section, we construct some simple relative trisections and relative trisection diagrams. The described relative trisection diagrams are pictured in Figure \ref{fig:smallrelative}.

\begin{example}[Diagrams on punctured $\mathbb{RP}^2$]\label{ex:mobiusdiagram}
Consider the relative trisection diagram $\mathcal{D}:=(M,\emptyset,\emptyset,\emptyset)$, where $M$ denotes the Mobius band. This is the smallest (i.e.\ unique relative trisection diagram with $\Sigma\cong M$) non-orientable relative trisection diagram. Following the procedure of Proposition \ref{prop:reldiagrams}, the manifold described by $\mathcal{D}$ is $M\times D^2\cong B^3\ttimes S^1$.

Now let $M_n$ denote $M$ with $n\ge 0$ disjoint open disks removed. Every simple closed curve on $M_n$ that has an annular neighborhood is separating, so again we conclude there is a unique relative trisection diagram on $M_n$: $(M_n,\emptyset,\emptyset,\emptyset)$. This diagram describes a relative trisection of $M_n\times D^2\cong \natural_{n+1} B^3\ttimes S^1$.
\end{example}

\begin{example}[Diagrams on punctured Klein bottles]\label{ex:kleindiagram}
Let $K_n$ denote the Klein bottle with $n>0$ open disks removed. The only standard pairs of curve sets on $K_n$ are the empty pair and a pair of parallel curves each compressing $K_n$ to a planar surface. We conclude that there are two relative trisection diagrams on $K_n$: $(K_n,\emptyset,\emptyset,\emptyset)$ and $(K_n,c,c,c)$, where $c$ denotes a nonseparating curve on $K_1$ with $\omega_1([c])=0$.  Respectively, these are diagrams of $K_n\times D^2\cong\natural_{n+1}B^3\ttimes S^1$ and $(S^3\ttimes S^1\setminus\mathring{B}^4)\natural(\natural_{n-1} B^3\ttimes S^1$).
\end{example}

\begin{example}[Diagram of $D^2\ttimes\mathbb{RP}^2$]\label{ex:rp2diskbundle}
The smallest irreducible relative trisection diagram $\mathcal{D}=(\Sigma;\alpha,\beta,\gamma)$ is pictured in Figure \ref{fig:smallrelative} (rightmost). The surface $\Sigma$ is the non-orientable surface with one boundary component and $\chi(\Sigma)=-2$. Each of $\alpha,\beta,\gamma$ consist of one curve, with each pair consisting of dual curves. According to Algorithm \ref{example:getkirbyfromrel} (discussed in the next subsection), the manifold $X^4$ described by $\mathcal{D}$ is $D^2\ttimes\mathbb{RP}^2$. The open book induced on $\partial X\cong S^2\ttimes S^1$ has page a Mobius band (and therefore trivial monodromy). This is the same (up to homeomorphism) open book as the one induced on $\boundary(B^3\ttimes S^1)$ by $(M,\emptyset,\emptyset,\emptyset)$ as in Example \ref{ex:mobiusdiagram}.
\end{example}

\subsection{Algorithms by example}

In this section, we describe three common procedures involving relative trisections.

\begin{algorithm}[Getting a Kirby diagram from a relative trisection]\label{example:getkirbyfromrel}
Suppose that $(\Sigma;\alpha,\beta,\gamma)$ is a relative trisection diagram. Let $\chi:=\chi(\Sigma)$. Say that each of $\alpha,\beta,\gamma$ consists of $n$ curves, and that there are $k'_{\alpha\beta}$ parallel $\alpha,\beta$ curves when $\alpha$ and $\beta$ are standardized. Similarly define $k'_{\beta\gamma}$ and $k'_{\gamma\alpha}$.

We will obtain a Kirby diagram for the manifold $X^4$ described by $(\Sigma;\alpha,\beta,\gamma)$ containing a page of the open book on $\boundary X^4$ induced by the described relative trisection. This is essentially the reverse of the procedure in \cite{reldiagrams}. In the case that $\Sigma$ is orientable, this procedure is described by \cite{surfacecomplements}.

An example of this procedure is illustrated in Figure \ref{fig:getrelkirby}, and is likely much more helpful than the ensuing wall of text.

\begin{enumerate}
\item[Step 0.] Standardize $\alpha,\beta$.
    \item [Step 1.] Embed $\Sigma$ into $\#_{k'_{\alpha\beta}}S^2\times S^1$ so that:
    \begin{itemize}
        \item $\Sigma$ is a (possibly one-sided) Heegaard surface for $\#_{k'_{\alpha\beta}}S^2\times S^1$.
        \item Let $\pi:\boundary(\nu(\Sigma))\to\Sigma$ be the induced double covering, with $S_1$ and $S_2$ the fundamental domains of $\pi$ induced by $\Sigma$. Then up to reordering $S_1$ and $S_2$, it must be the case that each curve in $\alpha$ (projected to $S_1$) and $\beta$ (projected to $S_2$) bounds a disk in the complement of $\Sigma$.
    \end{itemize}
   \textbf{(In words, embed $(\Sigma;\alpha,\beta)$ as a standard surface in $\#_{k'_{\alpha\beta}}S^2\times S^1$.)}
    \item [Step 2.] Draw $1-2n-\chi$ disjoint, properly embedded arcs $\eta$ in $\Sigma$ that are disjoint from $\alpha\cup\beta$ so that $\Sigma_\alpha\setminus\nu(\eta)$ and $\Sigma_\beta\setminus\nu(\eta)$ are disks, where $\Sigma_\alpha$ and $\Sigma_\beta$ denote the result of compressing $\Sigma$ along $\alpha$ or $\beta$, respectively. (We usually call $\eta$ a {\emph{cut system}} for $\alpha$ and for $\beta$.)
    \item [Step 3.] Compress $\Sigma$ along disks bounded by the $\alpha$ curves.
    \textbf{(In words, $\Sigma$ compressed along $\alpha$ will be the page of the open book.)}
    \item [Step 4.] Let $C_1,\ldots, C_{1-2n-\chi}$ be the curves in $\boundary\nu(\Sigma)$ double-covering $\eta_1,\ldots,\eta_{1-2n-\chi}$. Draw a dot on each $C_i$. These curves will describe additional 1-handles in the final Kirby diagram. Fix an orientation of $\partial\Sigma$. If a neighborhood of $\eta_i$ is an orientation-preserving band attached to $\boundary\Sigma$, then $C_i$ will be a standard dotted circle (orientatable $1$-handle). Otherwise, $C_i$ will correspond to a non-orientable 1-handle.
    \textbf{(In words, parallel $\alpha,\beta$ curves and arcs in a cut system for $\alpha,\beta$ yield 1-handles.)}
    \item [Step 5.] Let $\gamma'$ be a collection of $(n-k'_{\beta\gamma})$ $\gamma$ curves that are linearly independent in $H_1(\Sigma)$ from $\beta$. Frame each arc of $\gamma'\setminus\eta$ according to the surface framing of $\Sigma\setminus\nu(\eta)$. The $\gamma'$ curves can be pushed off $\Sigma$ to yield attaching circles of $2$-handles. \textbf{(In words, $\gamma$ curves that are dual to $\beta$ curves yield $2$-handles.)}
    \item [Step 6.] Replace the dotted circles with the standard $1$-handle notation. When $C_i$ describes a non-orientable $1$-handle, this handle is obtained by deleting a slice disk for $C_i$ as usual for dotted circles, and then cutting and regluing the resulting $S^2$ by an orientation-reversing map. This map fixes a whole circle on the $S^2$, which arrange to agree with the intersection of $S^2$ with $\Sigma$ and $\gamma$. Then we can replace $C_i$ with a non-orientable $1$-handle as usual.
    \item [Step 7.]We obtain a Kirby diagram $\mathcal{K}$ in which the framed $\gamma$ curves are the 2-handle attaching circles. If $k'_{\gamma\alpha}=0$, then we are done.
    \item [Step 8.] Standardize $\gamma,\alpha$ and perform the corresponding handle slides in $\mathcal{K}$. Then for each $\gamma$ curve $c$ that is parallel to an $\alpha$ curve $a$, add a 3-handle to $\mathcal{K}$ along the planar surface cobounded by $c$ and $a$, capped off by a cocore of the 1-handle corresponding to $a$ and core of the 2-handle corresponding by $c$, if they exist, or else just capped off by disks bounded by $c$ and $a$.
    {\bf{(In words, $\gamma$ curves that are parallel to $\alpha$ curves yield $3$-handles.)}}
    \item [Step 9.] We are now done. The Kirby diagram $\K$ describes $X^4$ and contains $\Sigma$ as a page of the open book on $\partial X^4$ induced by the relative trisection described by $(\Sigma;\alpha,\beta)$. Simplify as desired.
\end{enumerate}
\end{algorithm}

\begin{figure}
    \centering
     \scalebox{.7}{
\begingroup%
  \makeatletter%
  \providecommand\color[2][]{%
    \errmessage{(Inkscape) Color is used for the text in Inkscape, but the package 'color.sty' is not loaded}%
    \renewcommand\color[2][]{}%
  }%
  \providecommand\transparent[1]{%
    \errmessage{(Inkscape) Transparency is used (non-zero) for the text in Inkscape, but the package 'transparent.sty' is not loaded}%
    \renewcommand\transparent[1]{}%
  }%
  \providecommand\rotatebox[2]{#2}%
  \newcommand*\fsize{\dimexpr\f@size pt\relax}%
  \newcommand*\lineheight[1]{\fontsize{\fsize}{#1\fsize}\selectfont}%
  \ifx\svgwidth\undefined%
    \setlength{\unitlength}{382.96592412bp}%
    \ifx\svgscale\undefined%
      \relax%
    \else%
      \setlength{\unitlength}{\unitlength * \real{\svgscale}}%
    \fi%
  \else%
    \setlength{\unitlength}{\svgwidth}%
  \fi%
  \global\let\svgwidth\undefined%
  \global\let\svgscale\undefined%
  \makeatother%
  \begin{picture}(1,0.60284683)%
    \lineheight{1}%
    \setlength\tabcolsep{0pt}%
    \put(0,0){\includegraphics[width=\unitlength,page=1]{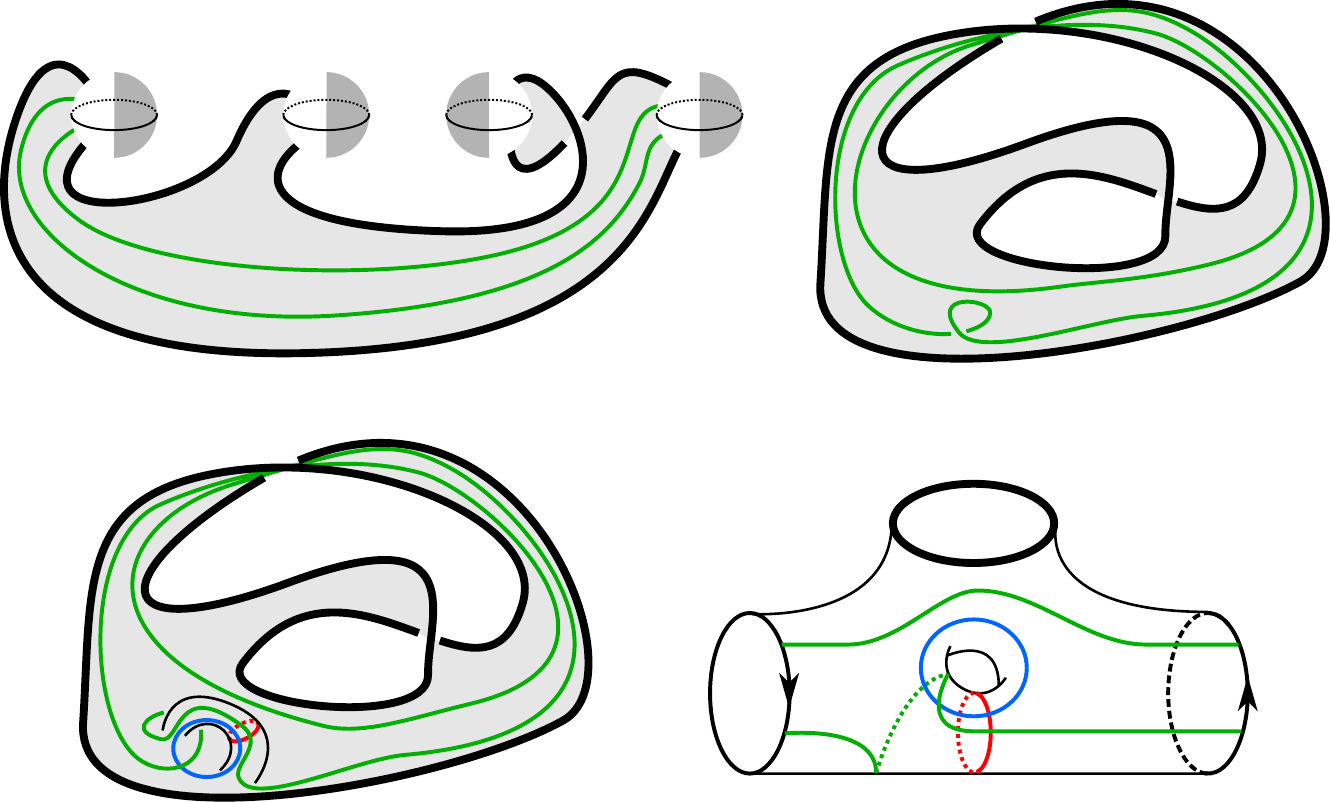}}%
    \put(0.23181866,0.50403791){\color[rgb]{0,0,0}\makebox(0,0)[lt]{\lineheight{1.25}\smash{\begin{tabular}[t]{l}A\end{tabular}}}}%
    \put(0.35413905,0.50403791){\color[rgb]{0,0,0}\makebox(0,0)[lt]{\lineheight{1.25}\smash{\begin{tabular}[t]{l}A\end{tabular}}}}%
    \put(0.07303292,0.50403791){\color[rgb]{0,0,0}\makebox(0,0)[lt]{\lineheight{1.25}\smash{\begin{tabular}[t]{l}B\end{tabular}}}}%
    \put(0.51341771,0.50403791){\color[rgb]{0,0,0}\makebox(0,0)[lt]{\lineheight{1.25}\smash{\begin{tabular}[t]{l}B\end{tabular}}}}%
    \put(0,0){\includegraphics[width=\unitlength,page=2]{getreltrisectionfromkirby.pdf}}%
    \put(0.12904227,0.36369374){\color[rgb]{0,0,0}\makebox(0,0)[lt]{\lineheight{1.25}\smash{\begin{tabular}[t]{l}$1$\end{tabular}}}}%
  \end{picture}%
\endgroup%
}
    \caption{As in Example \ref{example:getrelfromkirby}, we obtain a relative trisection diagram $(\Sigma;\alpha,\beta,\gamma)$ from a Kirby diagram $\K$ of $X^4$ containing a page $P$ of an open book $\mathcal{O}$ on $\partial X^4$. The diagram $(\Sigma;\alpha,\beta,\gamma)$ describes a relative trisection of $X^4$ that induces $\mathcal{O}$ on $\partial X^4$. We start in the top left with $\K$ and $P$. We have arranged that every band in some decomposition of $P$ runs over a distinct 1-handle of $\K$. Since there are no more $1$-handles, we do not obtain any parallel $\alpha,\beta$ curves. In the top right, we project the $2$-handle attaching circles to $P$ and do Reidemeister I moves so that the surface framing on each arc agrees with the handle framing. In the bottom left, we stabilize $P$ at each crossing of the projected circles to obtain an embedded circle $\gamma$. We add $\alpha,\beta$ curves at each stabilization. Since there are now an equal number of $\alpha,\beta,\gamma$ curves, we are done. We draw the whole diagram $(\Sigma,\alpha,\beta,\gamma)$ again in the bottom right.}
    \label{fig:getreltrisectionfromkirby}
\end{figure}

\begin{figure}
    \centering
     \scalebox{.85}{
\begingroup%
  \makeatletter%
  \providecommand\color[2][]{%
    \errmessage{(Inkscape) Color is used for the text in Inkscape, but the package 'color.sty' is not loaded}%
    \renewcommand\color[2][]{}%
  }%
  \providecommand\transparent[1]{%
    \errmessage{(Inkscape) Transparency is used (non-zero) for the text in Inkscape, but the package 'transparent.sty' is not loaded}%
    \renewcommand\transparent[1]{}%
  }%
  \providecommand\rotatebox[2]{#2}%
  \newcommand*\fsize{\dimexpr\f@size pt\relax}%
  \newcommand*\lineheight[1]{\fontsize{\fsize}{#1\fsize}\selectfont}%
  \ifx\svgwidth\undefined%
    \setlength{\unitlength}{453.43487116bp}%
    \ifx\svgscale\undefined%
      \relax%
    \else%
      \setlength{\unitlength}{\unitlength * \real{\svgscale}}%
    \fi%
  \else%
    \setlength{\unitlength}{\svgwidth}%
  \fi%
  \global\let\svgwidth\undefined%
  \global\let\svgscale\undefined%
  \makeatother%
  \begin{picture}(1,0.53165687)%
    \lineheight{1}%
    \setlength\tabcolsep{0pt}%
    \put(0,0){\includegraphics[width=\unitlength,page=1]{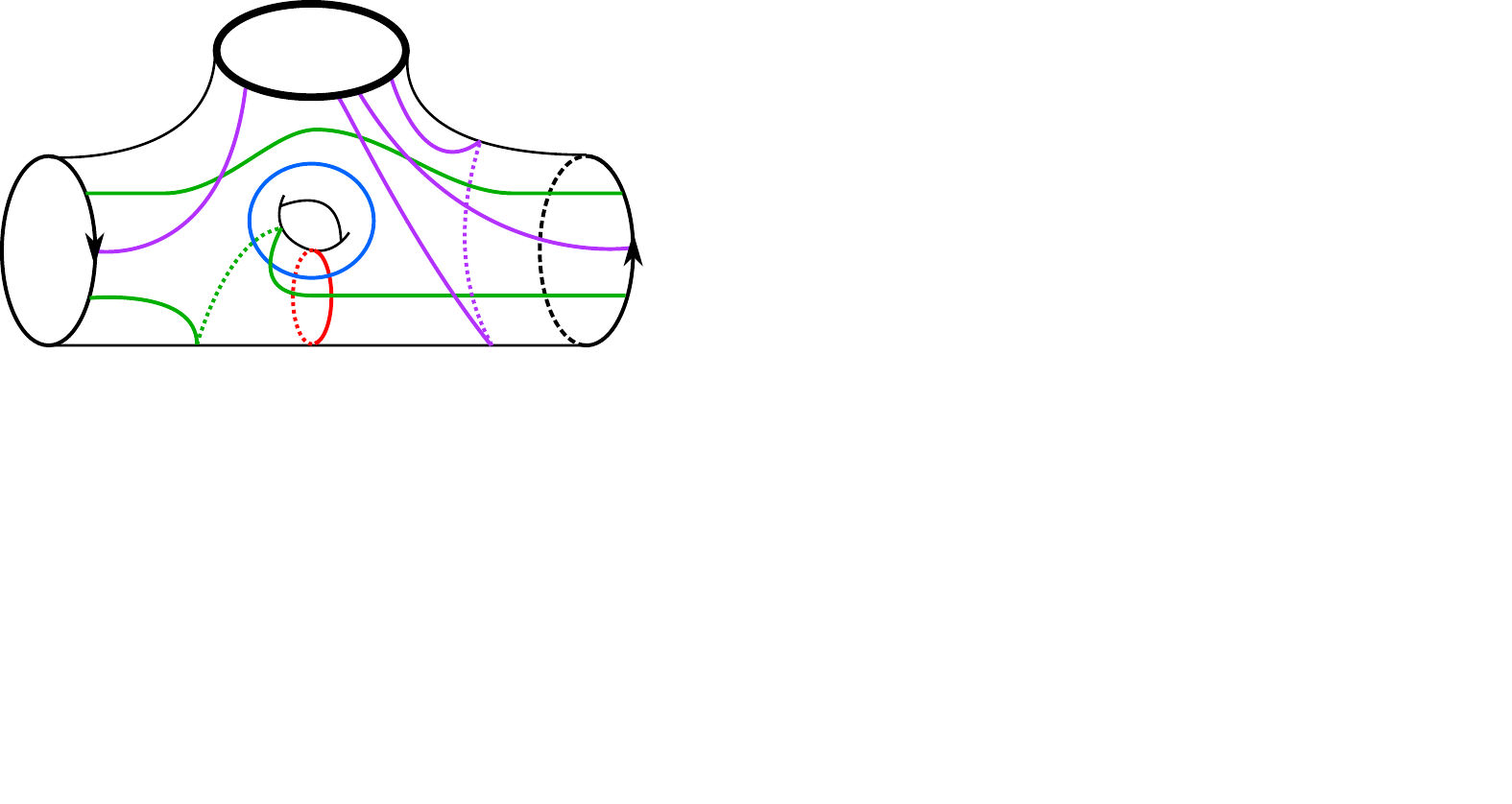}}%
    \put(0.27927709,0.47620784){\color[rgb]{0.50196078,0,0.50196078}\makebox(0,0)[lt]{\lineheight{1.25}\smash{\begin{tabular}[t]{l}$a=b$\end{tabular}}}}%
    \put(0,0){\includegraphics[width=\unitlength,page=2]{monodromyalgorithm.pdf}}%
    \put(0.84067733,0.47620784){\color[rgb]{0.50196078,0,0.50196078}\makebox(0,0)[lt]{\lineheight{1.25}\smash{\begin{tabular}[t]{l}$c$\end{tabular}}}}%
    \put(0,0){\includegraphics[width=\unitlength,page=3]{monodromyalgorithm.pdf}}%
    \put(0.27927706,0.17445518){\color[rgb]{0.50196078,0,0.50196078}\makebox(0,0)[lt]{\lineheight{1.25}\smash{\begin{tabular}[t]{l}$a'$\end{tabular}}}}%
    \put(0,0){\includegraphics[width=\unitlength,page=4]{monodromyalgorithm.pdf}}%
    \put(0.84425522,0.17445518){\color[rgb]{0.50196078,0,0.50196078}\makebox(0,0)[lt]{\lineheight{1.25}\smash{\begin{tabular}[t]{l}$a''$\end{tabular}}}}%
  \end{picture}%
\endgroup%
}
    \caption{We perform the monodromy algorithm as in Example \ref{example:monodromyalgorithm}. {\bf{Top left:}} We choose a cut system $a$ for $\alpha$ and $\beta$. {\bf{Top right:}} we slide $a$ over $\beta$ to be disjoint from $\gamma$ and call the resulting arcs $c$. {\bf{Bottom left:}} we slide $c$ over $\gamma$ to be disjoint from $\alpha$ and call the resulting arcs $a'$. The monodromy as an automorphism of $P_\alpha$ is now described by $\phi(a)=a'$. {\bf{Bottom right:}} we do further slides of $a'$ over $\alpha$ (i.e.\ isotopy in $P_{\alpha}$ to obtain arcs $a''$, which we see are isotopic rel boundary to $a$. We conclude that in this specific example, $\phi$ is isotopic to the identity.}
    \label{fig:monodromyalgorithm}
\end{figure}

Now we produce some relative trisections and diagrams thereof. The algorithm of \cite{casgaypin18} that produces a trisection from a Kirby diagram and page of an open book works verbatim, as illustrated in Figure \ref{fig:getreltrisectionfromkirby}. The monodromy algorithm of \cite{reldiagrams} also works directly. We illustrate this procedure in Figure \ref{fig:monodromyalgorithm}.

\begin{figure}
    \centering
     \scalebox{.6}{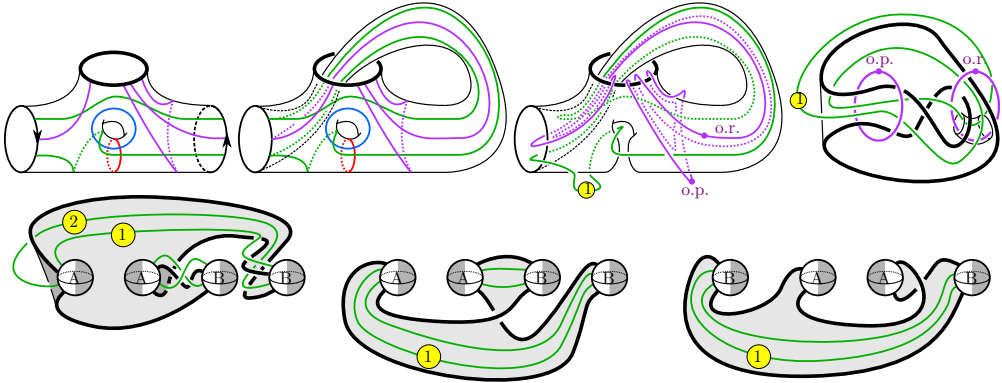}
    \caption{The procedure of obtaining a Kirby diagram from a relative trisection (see Example \ref{example:getkirbyfromrel}). In the top left, we have a relative trisection diagram $(\Sigma;\alpha,\beta,\gamma)$. We then embed the diagram in $S^3$ so that $(\Sigma;\alpha,\beta)$ is a standard one-sided surface. We draw arcs (in purple) on $\Sigma$ cutting $\Sigma_\alpha$ and $\Sigma_\beta$ into disks. Doubling these arcs yield 1-handles in a Kirby diagram (one of these 1-handles is not orientable). The unique $\gamma$ curve is homologically independent from $\beta$ and from $\alpha$, so it gives a 2-handle but not a 3-handle. We frame $\gamma$ according to the surface framing of $\Sigma$ to obtain the 2-handle attaching circle. Compressing $\Sigma$ along disks bounded by $\alpha$ yields a page of the open book on $\partial X^4$ induced by the relative trisection described by $(\Sigma;\alpha,\beta,\gamma)$. In the bottom row, we simplify this diagram further to see that the described 4-manifold is a $D^2\ttimes\mathbb{RP}^2\#B^3\times S^1$.}
    \label{fig:getrelkirby}
\end{figure}

\begin{algorithm}[Getting a relative trisection from a Kirby diagram]\label{example:getrelfromkirby}
Castro, Gay, and Pinz\'on-Caicedo \cite{casgaypin18} describe how to obtain a relative trisection of $X^4$ inducing a desired open book $\mathcal{O}$ on $\partial X^4$ from a Kirby diagram $\mathcal{K}$ of $X^4$ containing a page $P$ of $\mathcal{O}$. In the case that $X^4$ is not orientable, that procedure easily translates as follows (See Figure \ref{fig:getreltrisectionfromkirby}):
\begin{enumerate}
    \item[Step 0.] Fix a handle decomposition of $P$ with one $0$-handle and some number of bands. After adding cancelling $1$-, $2$-handle pairs of $\mathcal{K}$ and sliding $P$ over handles, we arrange that each band of $P$ goes over exactly one $1$-handle in $\mathcal{K}$, and at most one band of $P$ goes over any $1$-handle in $\mathcal{K}$.
    \item[Step 1.] For each $1$-handle of $\mathcal{K}$ that no band of $P$ runs over, fix a torus $T$ bounding a core of the $1$-handle in $\mathcal{K}$ and tube $P$ to $T$. Add parallel $\alpha,\beta$ curves on $T$ that are meridians of the core of the $1$-handle.
    \item[Step 2.] Project the $2$-handle attaching circles of $\mathcal{K}$ into $P$ (via the open book structure) as immersed curves. Break crossings as usual to indicate distance from $P$. Perform Reidemeister I moves to these projections as necessary so that the surface framing of $P$ on each attaching circle arc agrees with the handle framing, and so that each projected circle has at least one crossing (with itself or another circle).
    \item [Step 3.] At each crossing, attach an orientable tube to $P$ as in the bottom left of Figure \ref{fig:getreltrisectionfromkirby} to obtain a surface $\Sigma$. Add an $\alpha$ and $\beta$ dual pair on this tube. Have the overcrossing of the projected $2$-handle circle run over the tube and the undercrossing go between the feet of the tube, as pictured. The projected $2$-handle circles are now embedded; these are $\gamma$ curves.
    \item [Step 4.] If there are more $\alpha,\beta$ curves than $\gamma$ curves in the resulting diagram, add more $\gamma$ curves that are homologous to combinations of $\beta$ curves until there are an equal number of $\alpha,\beta,\gamma$ curves. (One can find such new $\gamma$ curves by banding together two $\beta$ curves that intersect a single $\gamma$ curve at adjacent points on $\gamma$.) Then $(\Sigma,\alpha,\beta,\gamma)$ is the desired relative trisection diagram.
\end{enumerate}
\end{algorithm}

\begin{algorithm}[Monodromy algorithm]\label{example:monodromyalgorithm}
In Figure \ref{fig:monodromyalgorithm}, we illustrate the monodromy algorithm of \cite{reldiagrams}, which allows one to compute the monodromy of the open book induced by a relative trisection diagram $(\Sigma;\alpha,\beta,\gamma)$. We compute the monodromy $\phi$ as an automorphism of the ``$\alpha$ page" $P_{\alpha}:= X_1\cap X_3\cap\partial X$ by identifying $P_\alpha$ with the result of compressing $\Sigma$ along $\alpha$.

\begin{enumerate}
    \item[Step 0.] Standardize $\alpha,\beta$. Let $a$ be a collection of disjoint properly embedded arcs in $\Sigma$, disjoint from $\alpha$ and from $\beta$, so that compressing $\Sigma\setminus(\nu(a))$ along $\alpha$ or along $\beta$ yields a disk. (We say that $a$ is a {\emph{cut system}} for $\alpha$ and $\beta$).
    \item[Step 1.] Do slides of $\beta,\gamma$ curves and slide $a$ over $\beta$ as necessary until transforming $a$ into arcs $c$ that are disjoint from $\gamma$. Note that the arcs $c$ might intersect $\alpha$ many times. If $\beta$ and $\gamma$ are standard then we can avoid slides of $\beta,\gamma$ curves, but otherwise we may have to perform many slides before obtaining $c$.
    \item[Step 2.] Let $\alpha'$ be another copy of $\alpha$. Do slides of $\gamma,\alpha'$ curves and slide $c$ over $\gamma$ as necessary until transforming $c$ into arcs $a'$ that are disjoint from $\alpha'$.
    \item[Step 3.] In practice, $\alpha'$ usually agrees with $\alpha$. However, we may have performed many slides. Now undo the slides to $\alpha'$ from the previous step to make $\alpha'$ again agree with $\alpha$ while simultaneously sliding $a'$ to remain disjoint from $\alpha'$. The monodromy $\phi$ is now described by $\phi(a)=a'$.
\end{enumerate}
\end{algorithm}

\begin{figure}
    \centering
    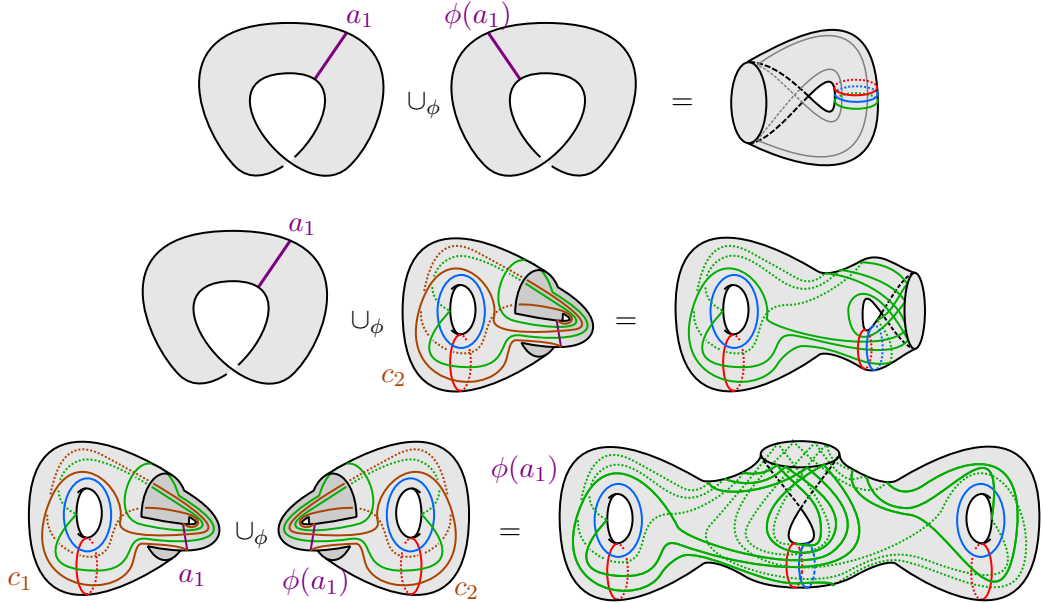
    \caption{As noted in Example \ref{ex:rp2diskbundle}, the pictured relative trisection diagrams of $B^3\ttimes S^1$ and $D^2\ttimes \mathbb{RP}^2$ induce equivalent open books on the respective boundaries. We can therefore glue the diagrams via the monodromy algorithm (Algorithm \ref{example:monodromyalgorithm}), as in Example \ref{ex:gluing}. {\bf{Top row:}} We glue relative trisections of $B^3\ttimes S^1$ to obtain a trisection of $S^3\ttimes S^1$ (as in Figure \ref{fig:basicexamples}(a)). {\bf{Middle row:}} We glue a relative trisection of $B^3\ttimes S^1$ to a relative trisection of $D^2\ttimes \mathbb{RP}^2$ to obtain a trisection of $\mathbb{RP}^4$ (as in Figure \ref{fig:basicexamples}(b)). {\bf{Bottom row:}} We glue relative trisections of $D^2\ttimes \mathbb{RP}^2$ to obtain a trisection of $S^2\ttimes \mathbb{RP}^2$ (as in Figure \ref{fig:twisteddouble}).}
    \label{fig:gluing}
\end{figure}

\begin{example}\label{ex:gluing}
The monodromy algorithm allows us to glue relative trisection diagrams. Let $\mathcal{D}_1=(\Sigma_1,\alpha_1,\beta_1,\gamma_1)$ and $\mathcal{D}_2=(\Sigma_2,\alpha_2,\beta_2,\gamma_2)$ be relative trisection diagrams of 4-manifolds $X,Y$ inducing open books $\mathcal{O}_1,\mathcal{O}_2$ on $\boundary X,\boundary Y$ respectively. Assume there exists a homeomorphism $\phi$ from $\boundary X$ to $\boundary Y$ taking the pages of $\mathcal{O}_1$ to the pages of $\mathcal{O}_2$. We say that $\mathcal{D}_1$ and $\mathcal{D}_2$ are {\emph{gluable}}. (The standard example is to fix a relative trisection diagram $\D_1$ of $X^4$ and then let $\D_2$ be obtained from $\D_1$ by an orientation-reversing automorphism. Then $\D_1$ and $\D_2$ can be glued to obtain a trisection diagram of the double of $X$ along its boundary.)

We glue the relative trisection diagrams by choosing a cut system of arcs $a_1$ in $\Sigma_1$ in the complement of $\alpha_1$, as in Algorithm \ref{example:monodromyalgorithm}. We perform the monodromy algorithm to obtain arcs $b_1$ disjoint from $\beta_1$ (if $\alpha,\beta$ are not already standardized, in which case we can choose $a_1$ so $a_1=b_1$ as we implicitly do in Algorithm \ref{example:monodromyalgorithm}) and $c_1$ disjoint from $\gamma$. Take the map $\phi:\boundary X\to\boundary Y$ to take the $\alpha$ page of the trisection induced by $\mathcal{D}_1$ to the $\alpha$ page of the trisection induced by $\mathcal{D}_2$. Then $a_2:=\phi(a_1)$ can be viewed as a cut system of arcs in $\Sigma_2$ for $\alpha_2$. We perform the monodromy algorithm again in $\Sigma_2$ to obtain arcs $b_2$ disjoint from $\beta_2$ and $c_2$ disjoint from $\gamma_2$.

Let $\Sigma:=\Sigma_1\cup_\boundary\Sigma_2$, using $\phi$ to identify boundary components of $\Sigma_1,\Sigma_2$. Then set \begin{align*}
\alpha:=\alpha_1\cup\alpha_2\cup(a_1\cup a_2),\\
\beta:=\beta_1\cup\beta_2\cup(b_1\cup b_2),\\
\gamma:=\gamma_1\cup\gamma_2\cup(c_1\cup c_2).\\
\end{align*}
Then $(\Sigma;\alpha,\beta,\gamma)$ is a trisection diagram for $X\cup_{\phi} Y$. We show various examples of this gluing operation in Figure \ref{fig:gluing}.
\end{example}

\begin{remark}
If $(\Sigma_1,\alpha_1,\beta_1,\gamma_1)$ and $(\Sigma_2,\alpha_2,\beta_2,\gamma_2)$ are gluable relative trisection diagrams, then $|\partial\Sigma_1|=|\partial\Sigma_2|$ and $\chi(\Sigma_1)\equiv\chi(\Sigma_2)\pmod{2}$. 
\end{remark}

\section{Bridge trisections of surfaces in non-orientable 4-manifolds}\label{sec:bridge}

In \cite{meizup17} and \cite{meizup18}, Meier and Zupan introduced \emph{bridge trisections} of knotted surfaces in 4-manifolds. This new way of studying knotted surfaces takes inspiration from a bridge splitting of a knot in a 3-manifold, and describes a knotted surface by three tangles in the 3-dimensional handlebodies of a trisection. In this section, we discuss how this theory can be modified to work in the non-orientable setting. While Meier--Zupan only consider orientable 4-manifolds, most of the theory carries over verbatim. For convenience, we include some short proofs, but focus on new examples. 

\subsection{Definitions}\label{sec:bridgetriexistence}

A collection $\D$ of properly embedded disks in a handlebody $Y$ is called \emph{trivial} if the disks in $\D$ are simultaneously boundary parallel. The following lemma is well-known in the orientable case, but we give a short proof for completeness. 

\begin{lemma}\label{lem:trivialdisks}
Suppose that $Y$ is a (possibly non-orientable) 4-dimensional handlebody of genus $p$, and that $U\subset \partial Y$ is an unlink. Then up to isotopy rel $L$, $L$ bounds a unique collection of boundary-parallel disks in $Y$. 
\end{lemma}

\begin{proof}[Proof]
Suppose that $\D$ and $\D'$ are two sets of disks such that $\partial \D=\partial \D'=L$, and let $\mathcal{S}\subset \partial Y$ be a collection of 2-spheres with the property that surgering $\partial Y$ along $\mathcal{S}$ yields $S^3$. Since $L$ is an unlink, we may isotope $\mathcal{S}$ so that $\mathcal{S}\cap L=\emptyset$. 

Given a 2-sphere $R$ in $\partial Y$, we could attach 3-and 4-handles to $\partial Y\times I$ to build a 4-manifold $Y'\cong Y$ with $\partial Y'$ and $\partial Y$ identified so that $S$ bounds a 3-ball into $Y'$. We conclude from \cite{laupoe72} or Theorem \ref{thm:NOLP} that there is a diffeomorphism rel boundary from $Y'$ to $Y$, so $R$ also bounds a 3-ball in $Y$. Using this fact, a standard innermost argument shows that we may assume that $\mathcal{S}\cap \D=\mathcal{S}\cap \D'=\emptyset$. Thus by cutting along $\mathcal{S}$, we reduce to the case that $Y\cong B^4$, and the result follows from \cite{liv82}. 
\end{proof}

\begin{definition}[\cite{meizup18}]\label{def:bridge}
Let $X^4$ be a closed, connected 4-manifold, and suppose that $X$ has a trisection $\T$ with sectors $X_1$, $X_2$, and $X_3$. Let $S\subset X$ be a smoothly embedded surface. A \emph{bridge trisection} of $S$ with respect to $\T$ is a decomposition $(X,S)=(X_1,\D_1)\cup (X_2,\D_2)\cup(X_3,\D_3)$, where $\D_i$ is a collection of trivial disks in $X_i$, and for $i\neq j$ the arcs $\tau_{ij}=\D_i\cap \D_j$ form a trivial (boundary-parallel) tangle in $H_{ij}$.
\end{definition}

Recall that a \emph{shadow} for a trivial arc $t$ in a handlebody $H$ is an embedded arc $s\subset\partial H$, such that $t$ and $s$ are isotopic in $H$ rel endpoints.  A collection of trivial arcs may admit non-isotopic sets of shadows (in $\partial H$), but these are related by slides of one shadow over another and over curves bounding disks into $H$. Bridge trisections may be described diagrammatically, via {\emph{shadow diagrams.}} 

\begin{definition}[\cite{meizup18}]
Let $X^4$ be a closed, connected 4-manifold, and suppose $S\subset X$ is a smoothly embedded surface. A $(g,k_i;b,c_i)$-{\emph{bridge trisection diagram}} (or {\emph{shadow diagram}}) for $S\subset X$ is a diagram $(\Sigma_g;\alpha,\beta,\gamma;s_\alpha,s_\beta,s_\gamma)$, where:
\begin{enumerate}
    \item $(\Sigma;\alpha,\beta,\gamma)$ is a $(g;k_i)$-trisection diagram with $2b$ additional marked points;
    \item The arcs $s_\alpha$,$s_\beta$, and $s_\gamma$ are three collections of $b$ shadows for trivial tangles $t_\alpha\subset H_\alpha$, $t_\beta\subset H_\beta$, and $t_\gamma\subset H_\gamma$, respectively;
    \item The pairwise unions $t_\alpha\cup t_\beta$, $t_\beta\cup t_\gamma$, and $t_\gamma\cup t_\alpha$ are $c_1$-, $c_3$-, and $c_2$- component unlinks, respectively.
\end{enumerate}
\end{definition}

These definitions extend verbatim to the non-orientable setting. By Lemma \ref{lem:trivialdisks}, any set of trivial disks for a bridge trisection of $K$ are isotopic, and so $K$ is completely determined by the arcs $t_\alpha$, $t_\beta$, $t_\gamma$. The proof of existence of bridge trisections is due to Meier and Zupan \cite{meizup18}, and the proof of uniqueness in the following sense is due in the standard trisection of $S^4$ to Meier and Zupan \cite{meizup17} and more generally to Hughes, Kim, and the first author \cite{hugkimmil20}. Both may be extended to the non-orientable setting with minor modifications. We will outline the key steps, but refer the reader to these references for more details. 

\begin{theorem}\label{thm:bridgeunique}
Let $X$ be a closed connected (but not necessarily orientable) 4-manifold with a trisection $\T$, and $S\subset X$ a smoothly embedded surface. Then $S$ may be isotoped to be in bridge trisected position with respect to $\T$. Moreover, if $S'$ is isotopic to $S$, then any two bridge trisections for $S$ and $S'$ become isotopic after some sequence of perturbations and deperturbations. 
\end{theorem}

The perturbation operation is defined in Section \ref{sec:perturbation}, and we comment on Theorem \ref{thm:bridgeunique} in Section \ref{sec:bridgeunique}, after discussing banded unlink diagrams.
\subsection{Existence and examples}

The proof of existence follows from the existence of banded unlink diagrams, which we briefly outline. 

\begin{definition}
Let $X$ be a closed connected 4-manifold, and $S\subset X$ a smoothly embedded surface. A function $f:X\to \mathbb{R}$ is a \emph{Morse function for the pair $(X,S)$} if $f$ is Morse, and $f\vert_S:S\to \mathbb{R}$ is also Morse. Moreover, $f$ is called \emph{self-indexing} if the image of all index $k$ critical points for $f$ are contained in $f^{-1}(k)$. We will always assume that $f$ has a unique index-0 critical points and a unique index-4 critical point.
\end{definition}

By a mild isotopy of $S$ in $X$, $f\vert_K$ may be assumed to be Morse, so we will always work with a Morse function of the pair $(X,S)$. We fix a gradient-like vector field $\nabla$ for $f$. For $M\subset X$, we will denote $M_{[a,b]}:=f^{-1}[a,b]\cap M$. The following proposition is well-known, and the proof carries over verbatim to the non-orientable setting. We give a short proof for convenience, since this gives an alternate proof of the existence of trisections for non-orientable 4-manifolds.

\begin{proposition}[\cite{gaykirby},\cite{meizup18}]\label{prop:handletri}
Suppose that $X$ is a closed and connected 4-manifold, and that $f:X\to\mathbb{R}$ is a self-indexing Morse function. Then $f$ induces a trisection of $X$. 
\end{proposition}

\begin{proof}[Proof sketch]
Suppose that $f$ has $n_i$ critical points of index $i$. The descending spheres of the index-2 critical points define the attaching link $L\subset X_{\{3/2\}}\cong \#_{n_1}S^2\times S^1$ or $\#_{n_1}S^2\ttimes S^1$ for the $2$-handles of $X$. Moreover, we can find a (possibly non-orientable) genus $g$ Heegaard splitting $X_{\{3/2\}}=H\cup_\Sigma H'$ of this level set so that $L$ is a core of $H$, i.e.\  every component of $L$ is dual to a properly embedded disk in $H$ that does not intersect any other component of $L$. The following decomposition now defines a $(g;n_1,g-n_2,n_3)$-trisection for $X$:
\begin{enumerate}
    \item $X_1=X_{[0,3/2]}\cup H_{[3/2,2]}$;
    \item $X_2=H'_{[3/2,5/2]}$;
    \item $X_3=H_{[2,5/2]}\cup X_{[5/2,4]}$.
\end{enumerate}
If desired, the trisection can be made balanced by stabilizing appropriately. 
\end{proof}

By an ambient isotopy of $S$ away from the critical points of $f$, we may also assume that $f\vert_S$ is self-indexing, and that all index 1 critical points of $f\vert_S$ are contained in $X_{\{3/2\}}$. In fact, we can always isotope $S$ to be in \emph{banded unlink position} (e.g.\  see \cite{meizup18}, \cite{hugkimmil20}).

\begin{definition}
Let $X$ be a 4-manifold, and let $f:X\to\mathbb{R}$ be a self-indexing Morse function. Let $S\subset X$ be a smoothly embedded surface. We say that $S$ is in \emph{banded unlink position} if:
\begin{enumerate}
    \item $f(S)=[1/2,5/2]$;
    \item $S_{\{1/2\}}\cap X_{\{1/2\}}$ and $S_{\{5/2\}}\cap X_{\{5/2\}}$ are collections of disjointly embedded disks;
    \item $S$ is \emph{vertical} on $(1/2,3/2)$ and $(3/2,5/2)$, i.e.\  $S\cap X_{(1/2,3/2)}=S_{(1/2,3/2)}$ and $S\cap X_{(3/2,5/2)}=S_{(3/2,5/2)}$;
    \item $S_{3/2}$ is a banded unlink, $(L,b)$, disjoint from the descending spheres of the index 2 critical points of $f$.
\end{enumerate}
\end{definition}

A banded unlink is an unlink $L\subset X_{\{3/2\}}$ together with attached \emph{bands} $b=\{b_1,...,b_m\}$ (corresponding to the index 1 critical points for $S$), with the property that resolving $L$ along $b$ is another unlink $L_b\subset X_{\{3/2\}}$. (Recall that a band attached to a link $L$ is a copy of $[0,1]\times[-\epsilon,\epsilon]$ meeting $L$ along $\{0,1\}\times[-\epsilon,\epsilon]$.) A \emph{banded unlink diagram} is a Kirby diagram $\K$ for $X$, together with $L$ and $b_1,\dots,b_m$. By Lemma \ref{lem:trivialdisks}, the triple $(\K,L,b)$ describes $S\subset X$. In \cite{hugkimmil20}, the authors study banded unlink diagrams in detail, and give a complete set of moves relating such diagrams. 

\begin{example}\label{ex:BUD}
Two simple examples of banded unlink diagrams are illustrated in Figure \ref{fig:bandedlinkdiagrams}. They describe familiar surfaces in the Kirby diagrams from Figure \ref{fig:fiberstructure}. In each case, the unknot bounds a disk in the $0$-handle of the Kirby diagram, and resolving each unlink with the band produces another unknot which bounds a disk in the $3$- and $4$-handle.
\begin{figure}[ht]
    \centering
    \includegraphics[width=0.7\textwidth]{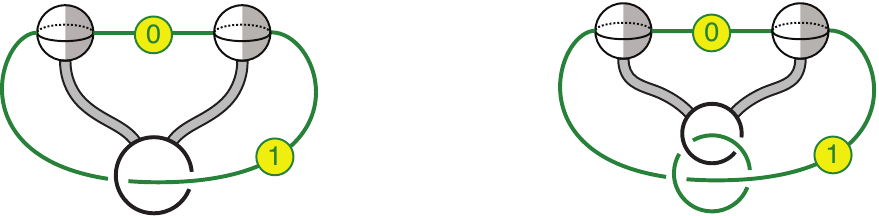}
    \caption{Two simple banded unlink diagrams. \textbf{Left}: $\mathbb{RP}^2\subset\mathbb{RP}^4$. \textbf{Right}: $\{\text{pt}\}\times\mathbb{RP}^2\subset  S^2\ttimes \mathbb{RP}^2$.}%
    \label{fig:bandedlinkdiagrams}
\end{figure}
\end{example}

Given a surface $S\subset X$ in banded unlink position, we may further modify it so that it intersects the trisection $\T$ from Proposition \ref{prop:handletri} in bridge position as in Definition \ref{def:bridge}. A banded unlink diagram of $S\subset X$ is in \emph{bridge position} with respect to the trisection $\T$ if:

\begin{enumerate}
    \item $L\subset X_{\{3/2\}}$ is in bridge position with respect to $\Sigma\subset X_{\{3/2\}}$;
    \item Each band is surface-framed with respect to $\Sigma$, i.e.\  for each $i$, $b_i\cap \Sigma$ is a single arc;
    \item The bands are \emph{dual} to a set of shadows $s_\alpha$ for $L\cap H$, i.e.\  $s_\alpha$ and $b$ intersect only at their ends, and $s_\alpha\cup (\Sigma\cap b_i)$ has no closed components. 
\end{enumerate}

\begin{theorem}
Suppose that $X$ is a 4-manifold and $f:X\to \mathbb{R}$ is a self-indexing Morse function. If $S\subset X$ is a smoothly embedded surface, then $S$ may be isotoped to lie in bridge position with respect to the trisection $\T$ induced by $f$.
\end{theorem}

Meier and Zupan only considered orientable manifolds, but their proof carries over verbatim for non-orientable manifolds, and so we will not repeat it here. If $S\subset X$ is in bridge position as above, then pushing the bands $b_i$ into $H$ produces a bridge trisection of $S$ with respect to $\T$ \cite[Lemma 3.1]{meizup17}. In other words, $\D_i=S\cap X_i$ is a trivial disk system. It is easy to see that $\D_1$ and $\D_3$ are trivial disk systems, and the duality condition guarantees that $\D_2$ is also a trivial disk system. Similarly, $t_\alpha$ and $t_\beta$ are trivial tangles by construction, and the duality condition guarantees that $t_\gamma$ is also a trivial tangle.

\begin{lemma}\label{lem:BUDfromtri}
Suppose that $X$ is a 4-manifold and $S\subset X$ is a smoothly embedded surface described by a $(g,k_i;b,c_i)$-bridge trisection diagram $(\Sigma;\alpha,\beta,\gamma;s_\alpha,s_\beta,s_\gamma)$. Then we may extract a banded unlink diagram $(\K,L,b)$ for $S\subset X$ via the following procedure:
\begin{enumerate}
    \item Standardize the $\alpha$ and $\beta$ curves. Also, standardize the shadows so that $s_\alpha\cup s_\gamma$ is a collection of $c_2$ embedded closed curves on $\Sigma$. 
    \item Produce a Kirby diagram $\K$ from $(\Sigma;\alpha,\beta,\gamma)$ (e.g.\  see Section \ref{sec:nonorientabletrisections} or \cite{gaykirby}). Moreover, let $L$ be the link obtained by pushing the arcs $s_\alpha$ and $s_\beta$ into $H_\alpha$ and $H_\beta$, respectively.
    \item\label{addbands} Let $\{s_1',\dots,s_{c_2}'\}$ be a sub-collection of arcs from $s_\gamma$ (one from each component of $s_\alpha\cup s_\gamma$). Add $b-c_2$ surface-framed arcs $b=\{b_i\}$ (bands) to $L$ by pushing $s_\gamma-\{s_1',\dots,s_{c_2}'\}$ into $H_\alpha$.
\end{enumerate}
\end{lemma}

If $b=1$ (as in the case of a \emph{doubly pointed trisection diagram} \cite{meizup18}), then no bands are added in Step \ref{addbands}, and a banded unlink diagram for the resulting 2-sphere may be obtained by simply pushing $s_\alpha\cup s_\beta$ into $H_\alpha\cup H_\beta$. In this case, we need only draw the endpoints of the shadows, since there is a unique way to connect the two bridge points in the complement of each cut system. 

The proof of Lemma \ref{lem:BUDfromtri} is the same as \cite[Lemma 3.3]{meizup17}, with one exception. A key step (\cite[Proposition 5.1]{meizup18}) in the correspondence between banded unlink diagrams is the fact every bridge splitting of the $n$-component unlink in $\#_kS^2\times S^1$ is \emph{standard}, i.e.\  a perturbation of $n$ copies of the $1$-bridge splitting of the unlink in $S^3$ stabilized with $k$ copies of $S^2\times S^1$. Since this relies on a result of \cite{bacsch05} which is only stated for the orientable case, we give a proof of this fact. 

\begin{lemma}\label{lem:unlink}
Let $L$ be an $n$-component unlink in $M:=\#_k S^2\ttimes S^1$. Fix a Heegaard surface $F$ in $M$ and assume that $L$ is in bridge position with respect to $F$. Then $L$ can be deperturbed with respect to $F$ until each component of $L$ intersects $F$ in two points.
\end{lemma}

\begin{proof}
If $M\cong S^3$, then this follows from \cite{meizup18}. Assume that the claim holds whenever $k<K$ for some fixed $K>0$, and $M=\#_K S^2\times S^1$. Also assume the claim is true in $M$ for any unlink of fewer than $n$ components.

Note that $M\setminus\nu(L)\cong \#_KS^2\ttimes S^1\#_n(S^1\times D^2)$, and each component of $M\setminus(\nu(F)\cup\nu(L))$ is a compression body. 

By Haken's lemma for 3-manifolds with boundary (see \cite{bonahon} or \cite{reducing}, but we recommend Chapter II of \cite{jaco}\footnote{The proof is only for closed manifolds, but the relative case is similar. Orientability does not play a role in the proof.}), there is an essential 2-sphere $S$ in $M\setminus L$ that intersects $F$ in a connected simple closed curve. 
Surger $M$ along $S$ (i.e.\ delete $\nu(S)$ and reglue two 3-balls). There is an induced Heegaard splitting of the resulting manifold $M'$.

If $M'$ is connected, then it is homeomorphic to $\#_{K-1} S^2\times S^1$ or $\#_{K-1} S^2\ttimes S^1$. In the first case, the claim holds by \cite{meizup18}; in the second it holds by inductive hypothesis.

If $M'=M'_1\sqcup M'_2$ is disconnected, then let $L_i=M'_i\cap L$. Since $S$ is essential in $M\setminus L$, for each $i$ either:
\begin{itemize}
    \item $M'_i\cong\natural_k S^2\times S^1$ for $k<K$, 
    \item $M'_i\cong\natural_k S^2\ttimes S^1$ for $k<K$,
    \item $L_i$ has fewer than $n$ components.
\end{itemize}
      In any case, the claim holds in $M'_1$ and $M'_2$ by \cite{meizup18} or inductive hypothesis.
\end{proof}

\begin{example}\label{ex:BT}
We end this section with two examples, producing the bridge trisection diagrams in Figure \ref{fig:nonORbridgetri}. A bridge trisection diagram for $\mathbb{RP}^2\subset \mathbb{RP}^4$ can be obtained by starting with the banded unlink diagram in Example \ref{ex:BUD}, and reversing the procedure in Lemma \ref{lem:BUDfromtri}. A doubly pointed diagram for $S^2\subset S^2\ttimes \mathbb{RP}^2$ may be found in the same way. One can use Lemma \ref{lem:BUDfromtri} to check that both diagrams are as advertised.

\begin{figure}[ht]
    \centering
    \includegraphics[width=\textwidth]{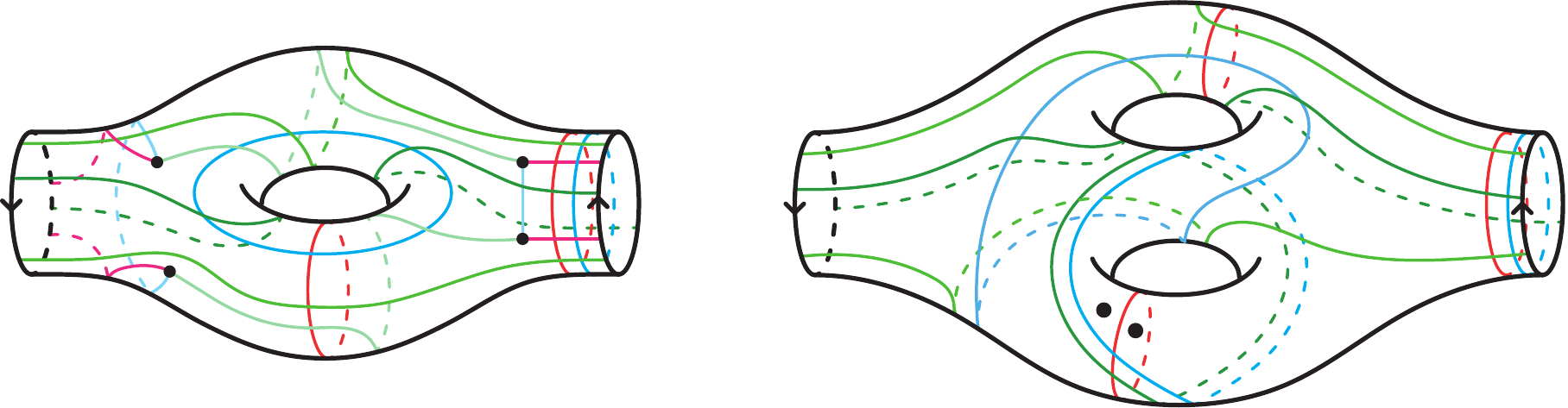}
    \caption{\textbf{Left:} A $(2,1;2,1)$- bridge trisection diagram for $\mathbb{RP}^2\subset \mathbb{RP}^4$. \textbf{Right:} A $(3,1;1,1)$- bridge trisection diagram for $S^2\times \{\text{pt}\}\subset S^2\times \mathbb{RP}^2$. Since there is a unique way to connect the two points in each of the cut systems, the shadows are not drawn. }%
    \label{fig:nonORbridgetri}
\end{figure}
\end{example}

\subsection{Perturbation}\label{sec:perturbation}

\begin{figure}
    \centering
    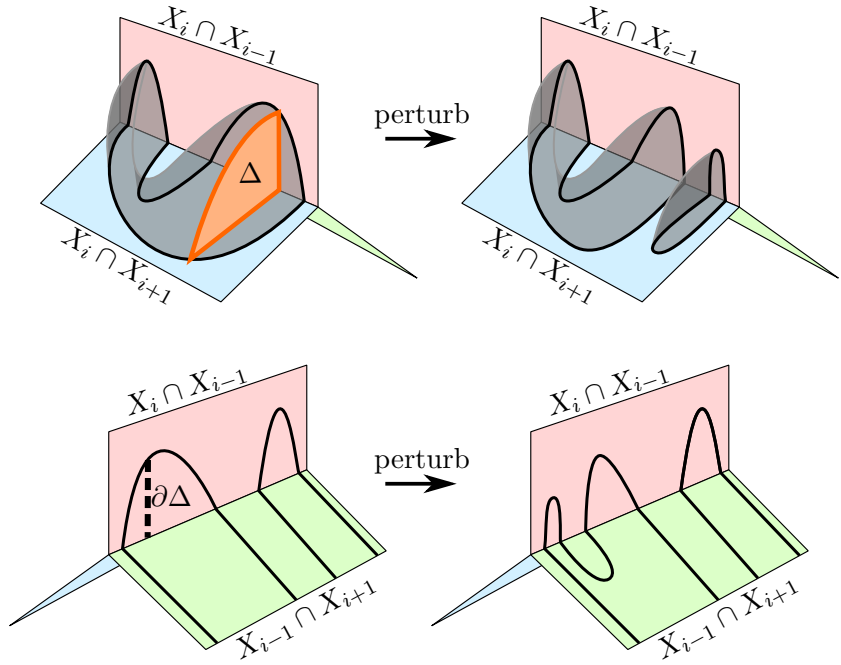
    \caption{{\bf{Left:}} part of a surface in bridge position, along with a disk $\Delta$ as in Definition \ref{def:perturb}. {\bf{Right:}} the result of perturbing the surface along $\Delta$.}
    \label{fig:perturbation}
\end{figure}

There is a natural way to perturb a surface in bridge position, analogous to perturbation of a knot in bridge position within a $3$-manifold. We illustrate this operation in Figure \ref{fig:perturbation}. 



\begin{definition}[\cite{meizup17,meizup18}]\label{def:perturb}
Let $S\subset X^4$ be a surface in bridge position with respect to a trisection $\mathcal{T}$ with sectors $X_1$, $X_2$, and $X_3$. Since $T_1=S\cap X_3\cap X_1$ and $T_2=S\cap X_1\cap X_2$ cobound boundary parallel disks in $X_1$, it follows that $T_1$ and $T_2$ have shadows $T'_1$ and $T'_2$ (respectively) in the central surface $\Sigma:=X_1\cap X_2\cap X_3=\boundary(X_3\cap X_1)=\boundary(X_1\cap X_2)$ that are disjoint in their interiors, and with $T'_1\cup T'_2$ an unlink bounding disks $D_1,\ldots, D_c$ in $\Sigma$. 
Then $S\cap X_1$ is isotopic rel $\boundary X_1$ to $D_1\sqcup\cdots\sqcup D_{c}$.

Let $\Delta$ be a disk in $X_1$ whose boundary can be decomposed into three arcs $\delta_1,\delta_2,\delta_3$ so that the following hold:
\begin{itemize}
    \item The arc $\delta_1$ is contained in $X_3\cap X_1$ with one endpoint on $S$ and the other on $\Sigma$. Moreover, projecting $\delta_1$ to $\Sigma$ yields an embedded arc with one endpoint on the interior of $T'_1$, one endpoint on $\Sigma\setminus S$, and the interior of $\delta_1$ disjoint from $T'_1$.
    \item The arc $\delta_2$ is contained in $S\cap X_1\cap X_2$ with one endpoint on $S$ and the other on $\Sigma$. Moreover, projecting $\delta_2$ to $\Sigma$ yields an embedded arc with one endpoint on the interior of $T'_2$, one endpoint on $\Sigma\setminus S$, and the interior of $\delta_2$ disjoint from $T'_2$,
    \item The arc $\delta_3$ is contained in $S$, and properly embedded in $X_1$.
\end{itemize}

Now let $S'$ be the surface obtained by compressing $S$ along $\Delta$. That is, frame $\Delta$ so that $(\delta_1\cup\delta_2)\times I\subset\boundary X_1$ and $\delta_3\times I\subset S$, and then let $S'$ be the result of performing a Whitney move on $S$ along $\Delta$; see Figure \ref{fig:perturbation}. We say that $S'$ is obtained from $S$ by {\emph{perturbation}}, and $S$ is obtained from $S'$ by {\emph{deperturbation}}.

We may exchange the roles of $X_1,X_2$, and $X_3$ cyclically when performing this operation, i.e.\ alternatively obtain $S'$ from this compressing operation in either $X_2$ or $X_3$. We still say $S'$ is obtained from $S$ by perturbation and that $S$ is obtained from $S'$ by deperturbation.
\end{definition}

The definition of perturbation may seem unwieldy, but it is formulated specifically to yield a surface in bridge position.

\begin{proposition}\cite[Lemma 5.2]{meizup17,meizup18}
Let $S$ be a surface in bridge position with respect to a trisection $\mathcal{T}$ with sectors $X_1$, $X_2$, and $X_3$. Let $S'$ be obtained from $S$ by perturbation, using a disk in $X_i$. Then $S'$ is in bridge position with respect to $\mathcal{T}$.
\end{proposition}

Perturbation of bridge trisections is conveniently very similar to perturbation of a banded link in bridge position (see Figure \ref{fig:bandperturb}). There is a correspondence between bridge trisections of $S$ and bridge-split banded unlink diagrams of $S$ \cite{meizup18,hugkimmil20}, in which the two notions of perturbations agree.

\subsection{Uniqueness}\label{sec:bridgeunique}

\begin{figure}
    \centering
      \scalebox{.9}{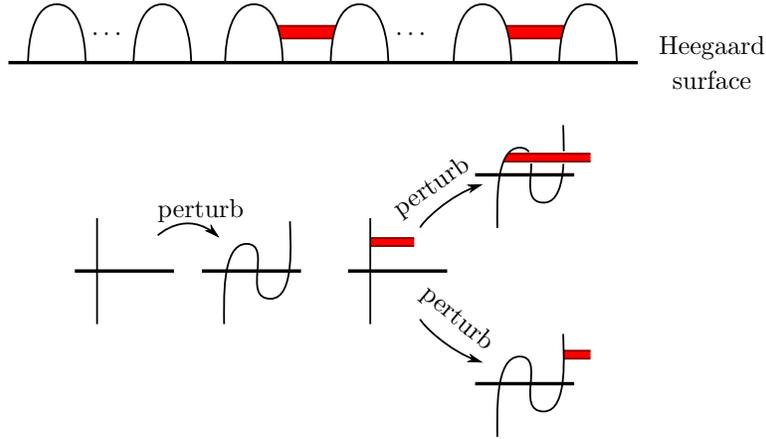}
    \caption{{\bf{Top row:}} A banded link $(L,b)$ in a Heegaard split 3-manifold $H_1\cup H_2$ is in {\emph{bridge position}} with respect to $H_1\cap H_2$ if each $(L,b)\cap H_i$ is isotopic to this picture. {\bf{Bottom:}} We illustrate how to perturb a banded link in bridge position. The inverse operation (when it yields a banded link in bridge position) is called deperturbation.}
    \label{fig:bandperturb}
\end{figure}

The proof of Theorem \ref{thm:bridgeunique} in \cite{hugkimmil20} almost works verbatim for non-orientable trisection. In that paper (as in \cite{meizup17}), the authors use the correspondence between bridge trisections and banded unlink diagrams, explicitly showing how to induce any band move on a banded unlink via perturbations and deperturbations on the corresponding bridge trisection. However, the final step of the proof, which shows that two bridge trisections corresponding to isotopic banded unlinks are equivalent after perturbations (see Figure \ref{fig:bandperturb}), requires the following additional lemma.

\begin{lemma}\label{uniquenesslemma}

Let $M$ be a non-orientable handlebody of genus $g$. Write $M=H_1\cup H_2$, where $H_1\cong\Sigma\times I$ is a collar neighborhood of $\boundary M$ and $H_2=\overline{M\setminus H_1}$.

Let $T$ be a boundary-parallel tangle in $M$ that is in bridge position with respect to $H_1\cap H_2$, (i.e.\ the closure of any component of $T\setminus (H_1\cap H_2)$ is an arc properly embedded in $H_1$ or $H_2$ which is either parallel to an arc in $H_1\cap H_2$ or of the form $\{\pt\}\times I$ in $H_1=\Sigma\times I$).

Then after a finite sequence of deperturbations applied to $T_1$ and an isotopy fixing $H_1\cap H_2$ setwise, $T$ can be taken to intersect $H_1$ only in arcs of the form $\{\pt\}\times I$.

\end{lemma}

When $M$ is an orientable handlebody, \ref{uniquenesslemma} is a theorem of Hayashi and Shimokawa \cite{hayashi} (see also \cite{alexthesis}). Their proof involves reducing the general case to when $M$ is a 3-ball. We repeat that argument for non-orientable handlebodies.

\begin{proof}[Proof of Lemma \ref{uniquenesslemma}]
If $M$ is the 3-ball, the claim follows from \cite{hayashi}. Assume that the claim holds whenever the genus of $M$ is less than $g$. 

Observe that the claim is trivially true when $T$ is the empty tangle. Assume that the claim holds whenever $T$ has fewer than $n$ components for some fixed $n>0$, and then take $T$ to be an $n$-component tangle.

Since $T$ is boundary-parallel, $M\setminus\nu(T)$ is a non-orientable handlebody. By Haken's lemma for 3-manifolds with boundary (as in Lemma \ref{lem:unlink}), there exists a disk $D$ neatly embedded in $M$ so that:
\begin{itemize}
    \item $D$ is disjoint from $T$,
    \item $D$ intersects $H_1\cap H_2$ in a connected simple closed curve,
    \item $D$ is separating and each component of $M\setminus\nu(D)$ is either not a 3-ball or contains at least one component of $T$.
\end{itemize}
Let $M_1$ and $M_2$ be the components of $M\setminus\nu(D)$, and $T_i=M_i\cap T$. Note that $M_i$ is a handlebody. Let $F_i$ be a surface in $M_i$ parallel to the boundary of $M_i$, with $F_i$ agreeing with $H_1\cap H_2$ away from $D$. Then $T_i$ is in bridge position in $M_i$ with respect to $F_i$. Moreover, for each $i$, either the genus of $M_i$ is less than $g$ (and $M_i$ may be orientable) or $T_i$ is a tangle of fewer than $n$ components. Inductively, $T_i$ can be deperturbed relative to $F_i$ (and hence relative to $H_1\cap H_2$ since $F_i$ agrees with $H_1\cap H_2$ near $T_i$), to intersect $H_1$ only in arcs of the form $\{\pt\}\times I$ as desired.

\end{proof}

\begin{proof}[Sketch of the proof of Theorem \ref{thm:bridgeunique}]
Let $S_1$ and $S_2$ be isotopic surfaces that are each in bridge position with respect to a trisection $\mathcal{T}$ of a closed non-orientable 4-manifold $X^4$. As mentioned at the end of Section \ref{sec:perturbation}, there is a correspondence between bridge trisections with respect to $\mathcal{T}$ and banded unlink diagrams in bridge position in a Heegaard-split Kirby diagram $\mathcal{K}$ related to $\mathcal{T}$. By the proof of \cite{hugkimmil20}, two banded unlink diagrams of a surface in a 4-manifold with respect to the same Kirby diagram are related by a sequence of band moves and isotopy. In \cite{meizup17} and \cite{hugkimmil20}, it is shown how to achieve band moves of the banded unlinks corresponding to $S_1$ and $S_2$ via a sequence of perturbations and deperturbations of $S_1$ and $S_2$. Then we can take the associated banded unlinks to be isotopic.

Now we describe an alternative perspective of the situation, in which we have isotoped the two banded unlinks to agree. We have a banded unlink $(L,b)$ in $\# S^2\ttimes S^1$ and two Heegaard splittings $H_1\cup_F H_2=H_1'\cup_{F'}H_2'$ of $\# S^2\ttimes S^1$ so that $(L,b)$ is in bridge position with respect to both splittings. 
We want to prove that $F$ and $F'$ become isotopic as bridge surfaces for $(L,b)$ after a finite sequence of perturbations; we will use the argument of \cite[Theorem 2.2]{alexthesis}. 
Let $C_1$ and $C_2'$ be wedges of circles so that $H_1,H_2'$ deformation retract on $C_1,C_2'$, respectively. Generically, we take $C_1$ and $C_2'$ to be disjoint. By isotoping $F$ to lie close to $C_1$ and $F'$ to lie close to $C_2'$, we may therefore take $F$ and $F'$ to be disjoint, with $H_1$ and $H_2'$ disjoint. Let $W:=H_1'\setminus H_1=H_2\setminus H_2'$. Then $W$ is homeomorphic to $F\times I$. Let $F^*=F\times\{1/2\}\subset W$. By Lemma \ref{uniquenesslemma} applied to the splitting $(W,H_1)$ of $W\cup H_1$ and the tangle $L\cap (W\cup H_1)$, we find that $F^*$ is obtained from $F$ by perturbation. On the other hand, by Lemma \ref{uniquenesslemma} applied to the splitting $(W,H'_2)$ of $W\cup H'_2$ and the tangle $L\cap (W\cup H'_2)$, we find that $F^*$ is obtained from $F'$ by perturbation (without moving $b$, as in Figure \ref{fig:bandperturb}). This completes the proof.

\end{proof}


\bibliography{references}
\bibliographystyle{alpha}
\end{document}